\documentclass[11pt]{amsart}
\usepackage {amsmath,amssymb, amscd, amsthm, epsfig, color}
\input xy
\xyoption {all}

\newtheorem{theorem}{Theorem}
\newtheorem {lemma}{Lemma}
\newtheorem {corollary}{Corollary}
\newtheorem {proposition}{Proposition}

\theoremstyle{definition}

\theoremstyle {definition}
\newtheorem {definition}{Definition}

\def \pr {{\mathbb {P}^r}}
\def \fscheme {{\overline {M}_{0,n} (X, \beta)}}
\def \stack {{\overline {\mathcal {M}}_{0,n}(\pr, d)}}
\def \fstack {{\overline {\mathcal {M}}_{0,n}(X, \beta)}}

\def \c {{\mathbb C}}
\def \cs {{\mathbb C^{\star}}}
\def \t {{\mathbb T}}

\def \proof {{\bf Proof. }}

\def \cgath {{\overline {\mathcal M}^{H}_{\alpha}(\pr, d)}}
\def \ogath {{\mathcal M^{H}_{\alpha}(\pr, d)}}
\def \tgath {{\mathcal {\widetilde M}_{\alpha}^{H}(\pr, d)}}
\def \odgath {{{\mathcal M}_{\alpha}^{\Delta, H}(\pr, d)}}
\def \cdgath {{{\mathcal {\overline M}_{\alpha}^{\Delta, H}(\pr, d)}}}

\begin{document}

\title[Tautological Classes on Moduli Spaces]{Tautological Classes on the Moduli Spaces of Stable Maps to $\mathbb P^r$ via Torus Actions}
\author {Dragos Oprea}
\address {Massachusetts Institute of Technology,}
\address {77 Massachusetts Avenue, Cambridge, MA 02139.}
\email {oprea@math.mit.edu}
\date{}
\maketitle

In our previous paper \cite {O}, we introduced the tautological rings of the genus zero moduli spaces of stable maps to homogeneous spaces $X$. We showed that in the case of $SL$ flags, all rational cohomology classes on the stable map spaces are tautological using methods from Hodge theory. The purpose of this note is to indicate a localization proof, in the spirit of Gromov-Witten theory, when $X$ is a projective space. 

To set the stage, we recall the definition of the tautological rings. The moduli stacks $\overline {\mathcal M}_{0,S} (\pr, d)$ parametrize $S$-pointed, genus $0$, degree $d$ stable maps to $\pr$. We use the notation $\stack$ when the labeling set is $S=\{1, 2, \ldots, n\}$. These moduli spaces are connected by a complicated system of natural morphisms, which we enumerate below: \begin {itemize} 
\item \text {forgetful morphisms}: $\pi:\overline {\mathcal M}_{0,S} (\pr, d) \to {\overline {\mathcal M}}_{0, T}(\pr, d)$ \text {defined for} $T\subset S$. 
\item \text {gluing morphisms which produce maps with nodal domains}, $$gl:{\overline {\mathcal M}}_{0, S_1\cup \{\bullet\}}(\pr, d_1) {\times_{\pr}} {\overline {\mathcal M}}_{0, \{\star\} \cup S_2} (\pr, d_2) \rightarrow \overline {\mathcal M}_{0, S_1\cup S_2}(\pr, d_1+d_2).$$ 
\item {\text evaluation morphisms to the target space}, $ev_i:\overline {\mathcal M}_{0,S} (\pr, d) \to \pr$ for all $i\in S$.\end {itemize}

\begin {definition} The genus $0$ tautological rings $R^{\star}(\stack)$ are the smallest system of subrings of $A^{\star}(\stack)$ such that:  
\begin {itemize}
\item The system is closed under pushforwards by the gluing and forgetful morphisms. 
\item All monomials in the evaluation classes $ev_i^{\star} \alpha$ where $\alpha \in A^{\star} (\pr)$ are in the system. 
\end {itemize}
\end{definition}

The localization theorem in \cite {EG} can be used to show that the equivariant Chow rings of $\stack$ are tautological after {\it inverting} the torus characters. Our goal is to prove the following stronger result, still making use of the torus action.\\ 

{\bf Theorem. }{\em All rational Chow classes of $\stack$ are tautological.\/ }\\

While this result is not surprising, we hope that its {\it proof} could be of interest. Localization is a popular theme in Gromov-Witten theory, used extensively since the early papers on the subject. Nonetheless, our approach is novel in two ways. First, we make use of a non-generic torus action on $\pr$ which fixes one point $p$ and a hyperplane $H$: $$t\cdot [z_0:z_1:\ldots:z_r]=[z_0:tz_1:\ldots:tz_r].$$
Secondly, we completely determine the Bialynicki-Birula {\it plus} decomposition of the $\it {stack}$ of stable maps which describes the flow of maps under this action. In addition, we show that the decomposition is {\it filterable}. As a consequence, we build up the stack of stable maps by adding cells in a {\it well determined order}. This is the algebraic analogue of the Morse stratification, whose cells can be ordered by the levels of the critical sets. A filterable decomposition also gives a way of computing the Poincare polynomials of the moduli spaces of stable maps from those of the fixed loci. This method works quite well in low codimension, as we will demonstrate in a future paper.

Note that the Bialynicki-Birula decomposition has not been established for {\it general} smooth Deligne-Mumford {\it stacks} with a torus action. However, in our case, we succeed to explicitly write it down in the context of Gathmann-Li stacks \cite {G} \cite{Li}. These stacks compactify the locus of marked maps with fixed contact orders with the hyperplane $H$. Our approach for constructing the plus cells applies whenever we have an equivariant etale affine atlas. 

We now explain the main result. Decorated graphs $\Gamma$ will be used to bookkeep the fixed loci, henceforth denoted $\mathcal F_{\Gamma}$. Their vertices correspond to components or points of the domain mapped entirely to $p$ or $H$, and carry numbered legs for each of the markings, and degree labels. The edges, also decorated by degrees, correspond to the remaining components. We repackage the datum of a decorated graph $\Gamma$ into an explicit fibered product $\overline {\mathcal Y}_{\Gamma}$ of Kontsevich-Manin and Gathmann spaces in equation $\eqref{ygamma}$. The theorem stated above is a consequence of the following stronger result summarizing the properties of the plus decomposition:

\begin {theorem} \label {main} The stack $\overline {\mathcal M}_{0, n}(\pr, d)$ can be decomposed into disjoint locally closed substacks $\mathcal F_{\Gamma}^{+}$ (the ``plus" cells of maps ``flowing" into $\mathcal F_{\Gamma}$) such that:
\begin {itemize}
\item [(1)] The fixed loci $\mathcal F_{\Gamma}$ are substacks of $\mathcal F^{+}_{\Gamma}$. There are projection morphisms $\mathcal F^{+}_{\Gamma}\to \mathcal F_{\Gamma}$. On the level of coarse moduli schemes, we obtain the {\it plus} Bialynicki-Birula decomposition of the coarse moduli scheme of $\overline {\mathcal M}_{0, n}(\pr, d)$ for a suitable torus action.
\item [(2)] The decomposition is filterable. That is, there is a partial ordering of the graphs $\Gamma$ such that $${\overline {\mathcal F^{+}_{\Gamma}}}\subset \bigcup_{\Gamma'\leq \Gamma} \mathcal F^{+}_{\Gamma'}.$$
\item [(3)] The closures of $\mathcal F^{+}_{\Gamma}$ are images of the fibered products $\overline {\mathcal Y}_{\Gamma}$ of Kontsevich-Manin and Gathmann spaces under the tautological morphisms.
\item [(4)] The codimension of $\mathcal F^{+}_{\Gamma}$ can be explicitly computed from the graph $\Gamma$. If $\mathfrak u$ is the number of $H$-labeled degree $0$ vertices of total valency $1$, and $\mathfrak s$ is the number of $H$-labeled vertices which have positive degree or total valency at least $3$, then the codimension is $d+\mathfrak s-\mathfrak u$.
\item [(5)] The rational cohomology and rational Chow groups of $\stack$ are isomorphic.
\item [(6)] (There exists a collection of substacks $\xi$ which span the rational Chow groups of $\mathcal F_{\Gamma}$ and) there exist closed substacks $\overline {\xi^{+}}$ supported in $\overline {\mathcal F^{+}_{\Gamma}}$ (containing the locus of maps flowing into $\xi$) which span the rational Chow groups of $\overline {\mathcal M}_{0, n}(\pr, d)$. The stacks $\overline {\xi^{+}}$ are images of fibered products of Gathmann spaces and tautological substacks of the Kontsevich-Manin spaces to $H$. Inductively, they can be exhibited as images of fibered products of Gathmann and Kontsevich-Manin spaces of maps to various linear subspaces of $\pr$. 
\item [(7)] All rational Chow classes on $\overline {\mathcal M}_{0, n}(\pr, d)$ are tautological.\end {itemize}
\end {theorem}

This paper is organized as follows. The first section contains preliminary observations about localization on the moduli spaces of stable maps and about the Gathmann stacks. In the second section we construct the Bialynicki-Birula cells on a general smooth Deligne Mumford stack with an equivariant atlas. We establish the ``homology basis theorem" under a general filterability assumption. The third section contains the main part of the argument. There, we identify explicitly the torus decomposition for the stacks $\stack$, and show its filterability. Finally, the last section proves the main results.

We would like to thank Professor Gang Tian for encouragement, support and guidance, and Professor A. J. de Jong for several helpful conversations. 

{\bf Conventions.} All schemes and stacks are defined over $\mathbb C$. All stacks considered here are Deligne-Mumford. $\t$ stands for the one dimensional torus, which usually will be identified with $\cs$ via a fixed isomorphism. For schemes/stacks $X$ with a $\t$-action, $X^{\t}$ denotes the fixed locus.

\section {Preliminaries.}

In this section we collect several useful facts about the fixed loci of the torus action on the moduli spaces of stable maps. We also discuss the Gathmann compactification of the stack of maps with prescribed contact orders to a fixed hyperplane.

\subsection {Localization on the moduli spaces of stable maps.} The main theme of this paper is a description of the flow of stable maps under the torus action on $\stack$. This flow is obtained by translation of maps under the action on the target $\pr$. Traditionally, actions with isolated fixed points have been used. As it will become manifest in the next sections, it is better to consider the following action which in homogeneous coordinates is given by: \begin {equation}\label {action} t\cdot [z_0:z_1:\ldots: z_r]=[z_0: tz_1: \ldots: tz_r], \; \text {for } t\in \c^{\star}.\end {equation} There are two fixed sets: one of them is the isolated point $p=[1:0:\ldots: 0]$ and the other one is the hyperplane $H$ given by the equation $z_0=0$. 
We observe that \begin {equation}\label{flow} \text {if } z\in \pr - H \text { then } \lim_{t\to 0} t\cdot z = p. \end {equation}

The torus fixed stable maps $f:(C, x_1, \ldots, x_n) \to \pr$ are obtained as follows. The image of $f$ is an invariant curve in $\pr$, while the images of the marked points, contracted components, nodes and ramification points are invariant i.e. they map to $p$ or to $H$. The non-contracted components are either entirely contained in $H$, or otherwise they map to invariant curves in $\pr$ joining $p$ to a point $q_H$ in $H$. The restriction of the map $f$ to these latter components is totally ramified over $p$ and $q_H$. This requirement determines the map uniquely. To each fixed stable map we associate a tree $\Gamma$ such that:  
\begin {itemize} 
\item The edges correspond to the non-contracted components which are not contained in $H$. These edges are decorated with degrees.

\item The vertices of the tree correspond to the connected components of $f^{-1}(p) \cup f^{-1} (H)$. These vertices come with labels $p$ and $H$ such that adjacent vertices have distinct labels. Moreover, the vertices labeled $H$ also come with degree labels, corresponding to the degree of the stable map on the component mapped to $H$ (which is $0$ if these components are isolated points). 

\item $\Gamma$ has $n$ numbered legs coming from the marked points.
\end {itemize}
We introduce the following notation for the graph $\Gamma$. 
\begin {itemize}
\item By definition, the flags incident to a vertex are all incoming legs and half edges.\vskip.03in

\item Typically, {\it $v$ stands for a vertex labeled $p$}. We let $n(v)$ be its total valency (total number of incident flags). \vskip.03in

\item Typically, {\it $w$ stands for a vertex labeled $H$}. We let $n(w)$ be its total valency. The corresponding degree is $d_w$. \vskip.03in

\item The set of vertices is denoted $V(\Gamma)$. We write $V$ and $W$ for the number of vertices labeled $p$ and $H$ respectively. \vskip.03in

\item The set of edges is denoted $E(\Gamma)$. The degree of the edge $e$ is $d_e$. We write $E$ for the total number number of edges.\vskip.03in

\item For each vertex $v$, we write $\alpha_v$ for the collection of degrees of the incoming flags along with the datum of their distribution on the flags. We agree that the degrees of the legs are $0$. We use the notation $d_v=|\alpha_v|$ for the sum of the incoming degrees.\vskip.03in

\item A vertex $w$ labeled $H$ of degree $d_w=0$ is called {\it unstable} if $n(w)\leq 2$ and {\it very unstable} if $n(w)=1$. The unstable vertices correspond to the zero dimensional components of $f^{-1}(H)$, and have the following interpretation:\vskip.03in

\begin {itemize}

\item the very unstable vertices come from unmarked smooth points of the domain mapping to $H$;

\item the unstable vertices with one leg come from marked points of the domain mapping to $H$;

\item the unstable vertices with two incoming edges come from nodes of the domain mapping to $H$.

\end {itemize}\vskip.03in
The vertices $w$ labeled $H$ of positive degree or with $n(w)\geq 3$ are {\it stable}. \vskip.03in

\item Let $\mathfrak s$ be the number of stable vertices labeled $H$, and $\mathfrak u$ be the number of very unstable vertices.
\end {itemize}

The fixed locus corresponding to the decorated graph $\Gamma$ will be denoted by $\mathcal F_{\Gamma}$. It can be described as the image of a finite morphism: \begin {equation} \label {zeta} \zeta_{\Gamma}: \prod_{v \text { labeled } p} {\overline {\mathcal M}}_{0,n(v)}\times \prod_{w \text { labeled } H} {\overline {\mathcal M}}_{0, n(w)}(H, d_w) \to \stack\end {equation}
To get the fixed locus we need to factor out the action of a finite group $A_{\Gamma}$ of automorphisms, which is determined by the exact sequence below. The last term is the automorphism group of the decorated graph $\Gamma$: $$1\to \prod_{e\in E(\Gamma)} \mathbb Z/{d_e \mathbb Z}\to A_{\Gamma} \to \text {Aut}_{\Gamma}\to 1.$$ 
The map $\zeta_{\Gamma}$ can be described as follows. 
\begin {itemize}
\item For each vertex $v$ labeled $p$ pick a genus $0$, $n(v)$-marked stable curve $C_v$ . 
\item For each vertex $w$ labeled $H$ pick a genus $0$ stable map $f_w$ to $H$ of degree $d_w$ with $n(w)$ markings on the domain $C_w$.
\item When necessary, we need to interpret $C_v$ or $C_w$ as points. 
\end {itemize}

A fixed stable map $f$ with $n$ markings to $\pr$ is obtained as follows. 
\begin {itemize} \item The component $C_v$ will be mapped to $p$. The components $C_w$ will be mapped to $H$ with degree $d_w$ under the map $f_w$. 
\item We join any two curves $C_v$ and $C_w$ by a rational curve $C_e$ whenever there is an edge $e$ of the graph $\Gamma$ joining $v$ and $w$. We map $C_e$ to $\pr$ with degree $d_e$ such that the map is totally ramified over the special points. 
\item Finally, the marked points correspond to the legs of the graph $\Gamma$. 
\end {itemize}

\subsection {Gathmann's moduli spaces.} Gathmann's moduli spaces are an important ingredient of our localization proof. We briefly describe them below, referring the reader to \cite {G} \cite {GT} for the results quoted in this section. 

Let $\alpha=(\alpha_1, \ldots, \alpha_n)$ be a $n$ tuple of non-negative integers and write $l(\alpha)=n$. We will usually assume that: $$|\alpha|=\sum \alpha_i=d.$$ The substack $\overline {\mathcal M}_{\alpha}^{H}(\pr, d)$ of $\stack$ parametrizes stable maps $f:C \to \pr$ with markings $x_1, \ldots, x_n$ such that: \begin {itemize}\item $f(x_i) \in H \text { for all } i \text { such that } \alpha_i>0$;
\item $f^{\star} H-\sum_{i} \alpha_i x_i$ is effective.\end {itemize} Gathmann showed that this is an irreducible, reduced, proper substack of the expected codimension $|\alpha|=\sum_{i} \alpha_i$ of $\stack$. 

We will prove later that the Gathmann stacks define tautological classes on the moduli spaces $\stack$. To this end we will make use of the recursive structure of the Gathmann stacks explained in the equations ($\ref {gathmann}$) and ($\ref {gathmann1}$) below. We describe what happens if we increase the multiplicities. We let $e_j$ be the elementary $n$-tuple with $1$ in the $j^{\text {th}}$ position and $0$ elsewhere. Then, we have the following relation in $A_{\star}(\stack)$: \begin {equation}\label {gathmann} [{\overline {\mathcal M}}_{\alpha+e_j}^{H}(\pr, d)]=-(\alpha_j \psi_j + ev_j^{\star} H) \cdot [{\overline {\mathcal M}}_{\alpha}^{H}(\pr, d)] + [\mathcal D_{\alpha,j}(\pr, d)]
\end {equation}
The correction terms $\mathcal D_{\alpha,j}(\pr, d)$ come from the boundary of the Gathmann stacks. These boundary terms account for the stable maps $f$ with one ``internal" component $C_0$ mapped to $H$ with some degree $d_0$ and with some multiplicity conditions $\alpha^{0}$ at the marked points points of $f$ lying on $C_0$. Moreover, we require that the point $x_j$ lie on $C_0$. There are $r$ (union of) components attached to the internal component at $r$ points. For all $1\leq i\leq r$, the map has degree $d_i$ on the component $C_i$ and sends the intersection point with the internal component $C_0$ to $H$ with multiplicity $m^i$. In addition, there are multiplicity conditions $\alpha^{i}$ at the marked points of $f$ lying on $C_i$. We require that the $d_i$'s sum up to $d$ and that the $\alpha^{i}$'s form a partition of the $n$-tuple $\alpha$.

\begin {figure}
\begin {center}
\includegraphics*[scale=.6]{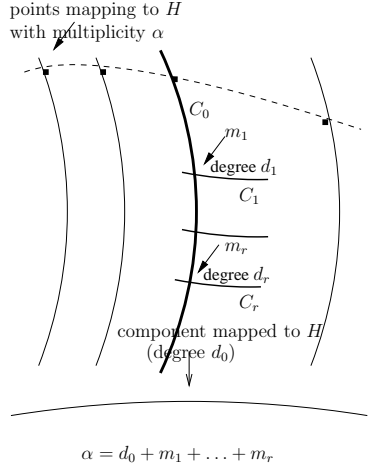}
\end {center}
\caption {A map in the boundary of the Gathmann compactification.}
\label {boundary}
\end {figure}
The boundary terms we just described are fibered products of lower dimensional Kontsevich-Manin and Gathmann stacks: $$\overline {\mathcal M}_{0, r+l(\alpha^0)} (H, d_{0})\times_{H^r}\prod_{i=1}^{r} \overline {\mathcal M}_{\alpha^i \cup m^i}^{H}(\pr, d_i).$$ Their multiplicities are found from the equation: \begin {equation}\label {gathmann1} [\mathcal D_{\alpha,j}(\pr, d)]=\sum \frac{m^1\ldots m^r}{r!}\left[\overline {\mathcal M}_{0, r+l(\alpha^0)} (H, d_{0})\times_{H^r}\prod_{i=1}^{r} \overline {\mathcal M}_{\alpha^i \cup m^i}^{H}(\pr, d_i)\right]\end {equation}

Jun Li extended Gathmann's construction to arbitrary genera and targets \cite {Li}. Li's stack $\mathfrak M^{H}_{0,\alpha}(\pr, d)$ parametrizes genus $0$ relative stable morphisms with contact order $\alpha$ to the pair $(\pr, H)$. Morphisms to level $k$ degenerations $\pr[k]$ of the target $\pr$ need to be considered. The collapsing maps $\pr[k]\to \pr$ induce a natural morphism $$\mathfrak M^{H}_{0,\alpha}(\pr, d)\to \mathcal {\overline M}_{0,n}(\pr, d)$$ whose image is the Gathmann stack $\overline {\mathcal M}_{\alpha}^H(\pr, d)$ \cite {GT}. Morphisms in the boundary of the Gathmann spaces which have components contained entirely in $H$ arise as collapsed images of morphisms to higher level degenerations of the target in Li's construction. An example is provided by the map shown in figure $\ref{boundary}$ which corresponds to a level $1$ relative morphism.

\section {The decomposition on smooth stacks with a torus action.}

In this section we will construct the Bialynicki-Birula cells of a smooth Deligne-Mumford stack with a torus action under the additional assumption that there exists an equivariant affine etale atlas. We show that the plus decomposition on the atlas descends to the stack. The existence of such an atlas should be a general fact, which we do not attempt to prove here since in the case of $\stack$ it can be constructed explicitly by hand. Finally, in lemma $\ref {homologybasis}$ we prove a ``homology basis theorem" for such stacks. Some (but not all) of the results presented in this section can in fact be proved from the corresponding statements for the coarse moduli schemes.

\subsection {The equivariant etale affine atlas.} In this subsection will construct an equivariant affine atlas for the moduli stack $\stack$. Fix an arbitrary $\t$-action on $\pr$ inducing an action by translation on $\stack$. Fix an identification $\t=\cs$. We may need to change this identification later. 

\begin {lemma} \label {atlas}
Possibly after lifting the action, there exists a smooth etale affine $\cs$-equivariant surjective atlas $S\to \stack$.  
\end {lemma}

{\proof} As a first step, we will find for any invariant stable map $f$, an equivariant etale atlas $S_f\to \stack$ whose image contains $f$. The construction in \cite {FP} shows that $\stack$ is a global quotient $[J/PGL(W)]$, thus giving a smooth surjective morphism $\pi:J\to \stack$. Here, $J$ is a quasiprojective scheme which is smooth since $\pi$ is smooth and $\stack$ is smooth. In fact, $J$ can be explicitly constructed as a locally closed subscheme of a product of Hilbert schemes on $\mathbb P(W)\times\pr$ for some vector space $W$. The starting point of the construction is an embedding of the stable map domain in $\mathbb P(W)\times \pr$. It is clear that the $\t$-action on the second factor equips $J$ with a $\t$-action such that the morphism $\pi:J\to \stack$ is equivariant. Moreover, from the explicit construction, it follows that $\pi^{\t}: J^{\t}\to \stack^{\t}$ is surjective. 

For any invariant $f$, there exists a $\t$-invariant point $j_f$ of $J$ whose image is $f$. It follows from \cite {S} that there exists an equivariant affine neighborhood $J_f$ of $j=j_f$ in $J$. The map on tangent spaces $d\pi: T_{j} J_f \to T_f \stack$ is equivariantly surjective. We can pick an equivariant subspace $V_f\hookrightarrow T_{j} J_f$ which maps isomorphically to $T_f\stack$. By theorem 2.1 in \cite {B1}, we can construct a smooth equivariant affine subvariety $S_f$ of $J_f$ containing $j$ such that $T_{j} S_f=V_f$. The map $\pi_f: S_f\to \stack$ is etale at $j$. Replacing $S_f$ to an equivariant open subset, we may assume $\pi_f$ is etale everywhere. Shrinking further, we can moreover assume that $S_f$ is equivariant smooth affine \cite {S}.

We consider the case of non-invariant maps $f$. We let $\alpha:\cs\to \stack$ be the equivariant {\it nonconstant} translation morphism: $$\cs\ni t\to f^{t}\in \stack.$$ Proposition 6 in \cite {FP} or corollary $\ref {eqfam}$ below show that, after possibly a base-change $\cs\to \cs$, we can extend this morphism across $0$. The image of $0\in \c$ under $\alpha$ is a $\t$-invariant map $F$ so we can utilize the atlas $S_F$ constructed above. We claim that the image of the atlas $\pi_F:S_F\to \stack$ contains $f$. Indeed, we consider the equivariant fiber product $C=\c\times_{\overline {\mathcal M}} S$. Since the morphism $C\to \stack$ is non-constant, the image of some closed point $j$ in $\stack$ is of the form $f^{t}$ for $t\neq 0$. Then, $f$ is the image of the closed point $t^{-1} j$. 

We obtained equivariant smooth affine atlases $S_f\to \stack$ whose images cover $\stack$. Only finitely many of them are necessary to cover $\stack$, and their disjoint union gives an affine smooth etale surjective atlas $S\to \stack$. 

\begin {corollary}\label {genatlas} Let $X$ be any convex smooth projective variety with a $\t$-action. There exists an equivariant smooth etale affine surjective atlas $S\to \fstack$ as in lemma $\ref {atlas}$.
\end {corollary}

{\proof} We embed $i:X\hookrightarrow \pr$ equivariantly, and base change the atlas $S$ constructed in the lemma under the closed immersion $i:\fstack \hookrightarrow \stack$. The convexity of $X$ is used to conclude that since $\fstack$ is smooth, the etale atlas $S$ is also smooth.

\subsection {The Bialynicki-Birula cells.} In this section we construct the Bialynicki-Birula cells for a smooth Deligne-Mumford stack $\mathcal M$ with a $\t$ action which admits an equivariant atlas as in proposition $\ref {atlas}$. This presupposes an identification $\t=\cs$. 

We will need the following observation. Let $f:X\to Y$ be an equivariant etale surjective morphism of smooth schemes (stacks) with torus actions. Let $Z$ be any component of the fixed locus of $Y$. Then $f^{-1}(Z)$ is union of components of $X^{\t}$ all mapping onto $Z$. Indeed, it suffices to show that the torus action on $f^{-1}(Z)$ is trivial. The $\t$-orbits in $f^{-1}(Z)$ need to be contracted by $f$ since $Z$ has a trivial action. Since the differential $df:TX\to TY$ is an isomorphism, it follows that all orbits are $0$ dimensional. They must be trivial since they are also reduced and irreducible. 

An affine fibration is a flat morphism $p:X\to Y$ which is etale locally trivial. Note that the transition functions need not be linear so $p$ is not necessarily a vector bundle. 

\begin {proposition}\label {stratification}
Let $\mathcal M$ be any smooth Deligne Mumford stack with a $\cs$-action and assume a $\cs$-equivariant affine etale surjective atlas $\pi:S\to \mathcal M$ has been constructed. Let $\mathcal F$ be the fixed substack and $\mathcal F_i$ be its connected components. Then $\mathcal M$ can be covered by locally closed disjoint substacks $\mathcal F_i^+$ which are affine fibrations over $\mathcal F_i$.
\end {proposition}

\proof Let $R=S\times_{\mathcal M} S$. The two etale surjective morphisms $s,t:R\to S$ together define a morphism $j: R\to S\times S$. It is clear that $R$ inherits a torus action such that $s,t$ are both equivariant. Moreover, since $\mathcal M$ is Deligne-Mumford, $j$ is quasi-finite, hence a composition of an open immersion and an affine morphism. Since $S$ is affine, it follows that $R$ is quasi-affine. As $s$ is etale, we obtain that $R$ is also smooth. 

If $F=S\times_{\mathcal M} \mathcal F$ then $F\hookrightarrow S$ is a closed immersion. Since $S\to \mathcal M$ is etale and equivariant, by the above observation, $F$ coincides with $S^{\t}$. Similarly $s^{-1}(F)$ and $t^{-1}(F)$ coincide with $R^{\t}$. Fixing $i$, we let $F_i=S\times_{\mathcal M}\mathcal F_i$. Then $F_i$ is union of components $F_{ij}$ of $S^{\t}$. Similarly, $s^{-1}(F_i)=t^{-1}(F_i)$ is a union of components $R_{ik}$ of $R^{\t}$. We will construct the substack $\mathcal F_i^{+}$ of $\mathcal M$ and the affine fibration $\alpha_i:\mathcal F_i^{+}\to \mathcal F_i$ on the atlas $S$. We will make use of the results of \cite {B1}, where a {\it plus} decomposition for quasi-affine schemes with a $\cs$-action was constructed. For each component $F_{ij}$, we consider its {\it plus} scheme $F_{ij}^{+}$; similarly for the $R_{ik}$'s we look at the cells $R_{ik}^{+}$. We claim that: $$s^{-1}(\cup_{j} F_{ij}^+)=t^{-1}(\cup_{j} F_{ij}^{+})=\cup_{k} R_{ik}^{+}$$ and we let $\mathcal F_i^{+}$ be the stack which $F_i^{+}=\cup_j F_{ij}^{+}\hookrightarrow S$ defines in $\mathcal M$. It suffices to show that if $R_{ik}$ is mapped to $F_{ij}$ under $s$, then $R_{ik}^{+}$ is a component of $s^{-1}(F_{ij}^{+})$. Let $r\in R_{ik}$. Then, using that $s$ is etale and the construction in \cite {B1}, we have the following equality of tangent spaces: $$T_r s^{-1}(F_{ij}^{+})=ds^{-1}\left(T_{s(r)}F_{ij}^{+}\right)=ds^{-1} \left((T_{s(r)}S)^{\geq 0}\right)=(T_r R)^{\geq 0}=T_r R^{+}_{ik}.$$ Here $V^{\geq 0}$ denotes the subspace of the equivariant vector space $V$ where the $\cs$-action has non-negative weights. The uniqueness result in corollary to theorem 2.2 in \cite {B1} finishes the proof under the observation that both $s^{-1}(F_{ij}^{+})$ and $R_{ik}^{+}$ are reduced. Note that the argument here shows that the codimension of $\mathcal F_{i}^{+}$ in $\mathcal M$ is given by the number of negative weights on the tangent bundle $T_r\mathcal M$ at a fixed point $r$.

To check that $\mathcal F_{i}^{+}\to \mathcal F_i$ is an affine fibration, we start with the observation that $F_{ij}^{+}\to F_{ij}$ are affine fibrations. We also need to check that the pullback fibrations under $s$ and $t$ are isomorphic: $$s^{\star}\left(\oplus_{j} F_{ij}^+\right)\simeq t^{\star}\left(\oplus_{j} F_{ij}^{+}\right)\simeq \oplus_{k} R_{ik}^{+}$$ The argument is identical to the one above, except that one needs to invoke corollary of proposition 3.1 in \cite {B1} to identify the fibration structure. Similarly, one checks the cocycle condition on triple overlaps.

Finally, we need to check that the trivializing open sets $U$ for $F_{ij}^{+}\to F_{ij}$ descend to the stack $\mathcal M$ i.e. we need to check that we can pick $U$ such that $s^{-1}(U)=t^{-1}(U)$. First, one runs the argument of lemma $2.2.3$ in \cite {AV}; we may assume that after replacing $U$ by an etale open, the groupoid $\xymatrix {R\ar[r]<.4ex>\ar[r]<-.4ex>&S}$ is (etale locally) given by $\xymatrix {U \times \Gamma\ar[r]<.4ex>\ar[r]<-.4ex> & U}$, where $\Gamma$ is a finite group acting on $U$. In this case, our claim is clear. 

\subsection {The homology ``basis" theorem.} In this subsection we will establish the ``homology basis theorem" (lemma $\ref {homologybasis}$) extending a result which is well known for smooth projective schemes \cite {C}. The proof does not contain any new ingredients, but we include it below, for completeness. We agree on the following conventions. The Chow groups we use are defined by Vistoli in \cite {V}, while the cohomology theory we use is defined for example in \cite {Be}. 

Let us consider a smooth Deligne Mumford stack $\mathcal M$ with a torus action whose fixed loci $\mathcal F_i$ are indexed by a finite set $I$, and whose Bialynicki-Birula cells $\mathcal F_i^{+}$ were defined above. We furthermore assume that the decomposition is filterable. That is, there is a partial (reflexive, transitive and anti-symmetric) ordering of the indices such that:
\begin {itemize}
\item [(a)] We have $\overline {\mathcal F_{i}^{+}}\subset \bigcup_{j\leq i} \mathcal F_{j}$;
\item [(b)] There is a unique maximal index $\mathfrak m\in I.$
\end {itemize}
Filterability of the Bialynicki-Birula decomposition was shown in \cite {B2} for projective schemes. For the stack $\stack$, filterability follows from the similar statement on the coarse moduli scheme. However, to prove the tautology of the Chow classes, we need the {\it stronger filterability condition} $(c)$, which we will demonstrate in the next section, and which does not follow from the known arguments:
\begin {itemize}
\item [(c)]  There is a family $\Xi$ of cycles supported on the fixed loci such that:
\begin {itemize} \item The cycles in $\Xi$ span the rational Chow groups of the fixed loci.
\item For all $\xi\in \Xi$ supported on $\mathcal F_i$, there is a {\it plus} substack $\xi^{+}=p_i^{-1}(\xi)$ (flowing into $\xi$) supported on $\mathcal F_{i}^{+}$; here $p_i:\mathcal F_{i}^{+}\to \mathcal F_i$ is the projection. We assume that $\xi^{+}$ is contained in a closed substack $\widehat {\xi^{+}}$ supported on $\overline {\mathcal F_i^{+}}$ (usually, but not necessarily, its closure) with the property: $$\widehat {\xi^{+}}\setminus \xi^{+} \subset \bigcup_{j < i} \mathcal F_{j}^+.$$ 
\end {itemize}
\end {itemize} 

\begin {lemma}\label {homologybasis} Assume that $\mathcal M$ is a smooth Deligne Mumford stack which satisfies the assumptions $(a)$ and $(b)$ above.
\begin {itemize}
\item [(i)] The Betti numbers $h^{m}(\mathcal M)$ of $\mathcal M$ can be computed as: 
$$h^{m}(\mathcal M) = \sum_{i} h^{m-2 n_i^{-}}(\mathcal F_i).$$ Here $n_i^{-}$ is the codimension of $\mathcal F_i^{+}$ which equals the number of negative weights on the tangent bundle of $\mathcal M$ at a fixed point in $\mathcal F_i$.
\item [(ii)] If the rational Chow rings and the rational cohomology of each fixed stack $\mathcal F_i$ are isomorphic, then the same is true about $\mathcal M$. 
\item [(iii)] Additionally, if assumption $(c)$ is satisfied, the cycles $\widehat {\xi^{+}}$ for $\xi \in \Xi$ span the rational Chow groups of $\mathcal M$. 
\end {itemize} 
\end {lemma}

{\proof } Thanks to item $(b)$, we can define an integer valued function $L(i)$ as the length of the shortest descending path from $\mathfrak m$ to $i$. Observe that $i<j$ implies $L(i)>L(j)$. Because of $(a)$, we observe that $$\mathcal Z_k=\bigcup_{L(i)> k} \mathcal F_i^{+}=\bigcup_{L(i)> k} \overline{\mathcal F_i^{+}}$$ is a closed substack of $\mathcal M$. Letting $\mathcal U_k$ denote its complement, we conclude that $$\mathcal U_{k-1}\hookrightarrow \mathcal U_k \text { and } \mathcal U_k \setminus \mathcal U_{k-1}\text { is union of cells } \bigcup_{L(i)=k}\mathcal F_{i}^{+}.$$ 

The Gysin sequence associated to the pair $(\mathcal U_k, \mathcal U_{k-1})$ is: $$\ldots \to \bigoplus_{L(i)=k} H^{m-2n^{-}_i}(\mathcal F_i^{+})=\bigoplus_{L(i)=k}H^{m-2n^{-}_i}(\mathcal F_i) \to H^{m}(\mathcal U_k) \to H^{m} (\mathcal U_{k-1})\to \ldots $$ One imitates the usual argument for smooth schemes in \cite {AB}, \cite {Ki} to prove that the long exact sequence splits. Item $(i)$ follows by estimating the dimensions. 

To prove $(ii)$, we compare all short exact Gysin sequences to the Chow exact sequences (for $m$ even) and use the five lemma:
\begin {equation}\label {chowgys}{
\xymatrix{
{}& \bigoplus_{L(i)=k}A^{m/2-n_i^{-}}(\mathcal F_i)\ar[r]\ar[d]& A^{m/2}(\mathcal U_k) \ar[r]\ar[d]& A^{m/2}(\mathcal U_{k-1})\ar[r]\ar[d]& 0\\
0\ar[r] & \bigoplus_{L(i)=k}H^{m-2n_i^{-}}(\mathcal F_i)\ar[r]& H^{m}(\mathcal U_k) \ar[r] & H^{m}(\mathcal U_{k-1})\ar[r]& 0}}
\end {equation} The first exact sequence exploits the surjectivity of the pullback $p_i^{\star}:A_{\star}(\mathcal F_i)\to A_{\star}(\mathcal F_{i}^{+})$ which follows by Noetherian induction from the local triviality of the affine bundle $p_i:\mathcal F_{i}^{+}\to \mathcal F_i$. The similar statement in cohomology is proved for instance in \cite {Be}. 

Finally, for $(iii)$, we use $\eqref {chowgys}$ to show inductively that $$\text { the cycles } \widehat {\xi^{+}}\cap \mathcal U_k \text { for } \xi \in \Xi \text { supported on } \mathcal F_i \text{ with } L(i)\leq k \text { span } A^{\star}(\mathcal U_k).$$ Condition $(c)$ is used to prove that the image in $\mathcal U_k$ of a cycle $\xi$ supported on a fixed locus $\mathcal F_i$ with $L(i)=k$ is among the claimed generators: $$\xi^{+}=\widehat {\xi^{+}}\cap \mathcal U_k.$$

\section {The Bialynicki-Birula decomposition on $\stack$.}
In the previous section we constructed the Bialynicki-Birula plus cells on the stack of stable maps $\stack$. In this section, we identify the decomposition explicitly. We start by analyzing the $\cs$-flow of individual stable maps. We will relate the decomposition to Gathmann's stacks in the next subsection. Finally, we will prove the filterability condition $(c)$ needed to apply lemma $\ref {homologybasis}$.

\subsection {The flow of individual maps.} To fix the notation, we let $f:(C, x_1, \ldots, x_n)\to \pr$ be a degree $d$ stable map to $\pr$. We look at the sequence of translated maps: $$f^t:(C, x_1, \ldots, x_n)\to \pr, \;\;f^t(z)=tf(z).$$ By the ``compactness theorem," this sequence will have a stable limit. We want to understand this limit $F=\lim_{t\to 0} f^t.$ 

To construct $F$ explicitly we need to lift the torus action $t\to t^{D}$ where $D=d!$. Henceforth, we will work with the lifted action: $$t\cdot [z_0:z_1:\ldots:z_r]=[z_0:t^{D}z_1:\ldots:t^{D}z_r].$$ We seek to construct a family of stable maps $G:\mathcal X\to \pr$ over $\c$, whose fiber over $t\neq 0$ is $f^t$ and whose central fiber $F:\mathcal C\to \pr$ will be explicitly described below.
\begin {equation}\label {family}
\xymatrix{
\mathcal C\ar@/^1pc/[rr]|F\ar[d]\ar[r]&\mathcal X\ar[r]^{G}\ar[d]^{\pi} & \pr \\
0\ar[r]\ar@/^1pc/[u]^{x_i}&\c\ar@/^1pc/[u]^{x_i}}
\end {equation}

First we assume that the domain $C$ is an irreducible curve. In case $f$ is mapped entirely to $H$, the family $\eqref {family}$ is trivial and $F=f$. 

Otherwise, $f$ intersects the hyperplane $H$ at isolated points, some of them possibly being among the marked points. We make a further simplifying assumption: we may assume that all points in $f^{-1}(H)$ are marked points of the domain. If this is not the case, we mark the remaining points in $f^{-1}(H)$ thus getting a new stable map $\overline f$ living in a moduli space with more markings $\overline {\mathcal M}_{0,n+k} (\pr, d)$. We will have constructed a family $\eqref {family}$ whose central fiber is $\overline F=\lim_{t\to 0} \overline f^t$.  A new family having $f^t$ as the $t$-fiber is obtained by forgetting the markings. We use a multiple of the line bundle $$\omega_{\pi}(\sum_{i} x_i)\otimes G^{\star} \mathcal O_{\pr}(3)$$ to contract the unstable components of the central fiber. Thus, we obtain the limit $F$ from $\overline F$ by forgetting the markings we added and stabilizing.

Henceforth we assume that all points in $f^{-1}(H)$ are among the markings of $f$, and $f$ is not a map to $H$. Let $s_1, \ldots, s_k$ be the markings which map to $H$, say with multiplicities $n_1, \ldots, n_k$ such that $\sum n_i = d$. We let $t_1, \ldots, t_l$ be the rest of the  markings. We let $q_i = f(s_i)$. The following lemma will be of crucial importance to us. The method of proof is an explicit stable reduction, and it is similar to that of proposition 2 in \cite {KP}.

\begin {lemma}\label {limit}
Let $F$ be the following stable map with reducible domain: 
\begin {itemize}
\item The domain has one component of degree $0$ mapped to $p$. This component contains markings $T_1, \ldots,  T_l$.
\item Additionally, there are $k$ components $C_1, \ldots, C_k$ attached to the degree $0$ component. The restriction of $F$ to $C_i$ has degree $n_i$, its image is the line joining $p$ to $q_i=f(s_i)$ and the map is totally ramified over $p$ and $q_i$. 
\item Moreover, if we let $S_i=F^{-1}(q_i)\cap C_i$, then $S_1, \ldots, S_k, T_1, \ldots, T_l$ are the marked points of the domain of $F$.
\end {itemize}
Then, the stabilization of $F$ is the limit $\lim_{t\to 0} f_t$. 
\end {lemma} 

{\bf Proof.} It suffices to exhibit a family as in $\eqref {family}$. We let $f^0, \ldots, f^r$ be the homogeneous components of the map $f$. We let $C$ be the domain curve with coordinates $[z:w]$. The assumption about the contact orders of $f$ with $H$ shows that $f^{0}$ vanishes at $s_1, \ldots, s_k$ of orders $n_1, \ldots, n_k$ with $\sum_{i} n_i = d$. 

There is a well defined map $G_0:\cs \times C {\to} \pr$ given by: $$(t, [z:w])\mapsto \left[f^0(z:w): t^{D} f^1(z:w):\ldots: t^{D} f^{r}(z:w)\right].$$ The projection map $\pi: \cs\times C \to \cs$ has constant sections $s_1, \ldots, s_k, t_1, \ldots, t_l$. It is clear that $G_0$ can be extended to a map $$G_0:\c\times C \setminus\bigcup_{i} \left(\{0\}\times \{s_i\}\right) \to \pr.$$ A suitable sequence of blowups of $\c \times C$ at the points $\{0\}\times \{s_i\}$ will give a family of stable maps $G:\mathcal X {\to} \pr$ as in $\eqref {family}$. 
\begin {figure}
\begin {center}
\includegraphics*[scale=.7]{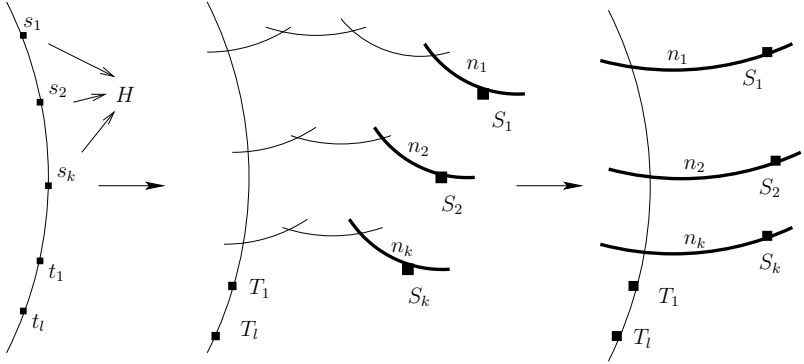}
\end {center}
\caption {Obtaining the stable limit.}
\label{evolution}
\end {figure}

It is useful to understand these blowups individually. It suffices to work locally in quasi-affine patches $U_i$ near $s_i$ and then glue. An affine change of coordinates will ensure $s_i=0$. For $n=n_i$, we write $f^0=z^n h$. We may assume that on $U_i$, $h$ does not vanish and that $f_1, \ldots, f_r$ do not all vanish. Let $D=n \cdot e$. We will perform $e$ blowups to resolve the map $G_0:U_i\times \mathbb P^{1}\setminus\{(0, [0:1])\}\to \pr$: $$(t, [z:w])\to \left [z^n h(z:w): t^{ne} f_1(z:w):\ldots: t^{ne} f_r(z:w)\right].$$ The blowup at $(0,[0:1])$ gives a map: $$\xymatrix{G_1: \mathcal X_1\ar@{.>}[r]& \pr}.$$ In coordinates, $$\mathcal X_1=\{(t, [z:w], [A_1, B_1]) \text { such that } A_1 z= t B_1 w \}$$ and $$G_1=\left[B_1^{n}h(tB_1:A_1): t^{ne -n} f_1(tB_1 : A_1):\ldots:  t^{ne-n} f_r(t B_1 : A_1)\right].$$
The map is still undefined at $t=0$ and $B_1=0$ so we will need to blow up again. After the $k^{\text {th}}$ blowup, we will have obtained a map: $$\xymatrix{G_k: \mathcal X_k\ar@{.>}[r]& \pr},$$ which in coordinates becomes: $$\mathcal X_k=\{(t, [z:w], [A_1, B_1], \ldots, [A_k:B_k])\;|\; A_1 z= t B_1 w,\; A_{i+1} B_i=t A_i B_{i+1}\}$$ 
$$G_k=\left[ B_k^{n} h(t^k B_k: A_k): t^{ne-nk} f_1(t^k B_k: A_k):\ldots :  t^{ne-nk} f_r(t^k B_k: A_k)\right].$$ After the $e^{\text {th}}$ blow up we obtain a well defined map. This map is constant on the first $e-1$ exceptional divisors (hence they are unstable). On the $e^{\text {th}}$ exceptional divisor the map is given by: $$G_e=\left[B_e^{n}h(0:1):A_e^{n} f_1(0:1):\ldots: A_e^{n} f_r(0:1)\right].$$ There, the map is totally ramified over two points in its image. It is easy to check that the sections $s_1, \ldots, s_k, t_1, \ldots, t_l$ extend over $t=0$ as claimed in the lemma. 

We obtain a family $G: \mathcal X \to \pr$ of maps parametrized by $\c$ as in $\eqref {family}$. The profile of the central fiber is the middle shape in figure $\ref {evolution}$. There are unstable components coming from the exceptional divisors which need to be contracted successively to obtain the final limit we announced. This completes the proof.

We consider the case when the domain curve is not irreducible. Assume that the stable map $f$ is obtained by gluing maps $f_1$ and $f_2$ with fewer irreducible components at markings $\star$ and $\bullet$ on their domains with $f_1(\star)=f_2(\bullet)$. Inductively, we will have constructed families $\eqref {family}$ of stable maps over $\c$ whose fibers over $t\neq 0$ are $f_1^t$ and $f_2^t$. We glue the two families together at the sections $\star$ and $\bullet$ thus obtaining a family whose fiber over $t$ is $f^{t}$. The argument above proves that the limit for reducible maps can be obtained by taking the limits of each irreducible component and gluing the limits together along the corresponding sections. 

\begin {figure}
\begin {center}
\includegraphics*[scale=.7]{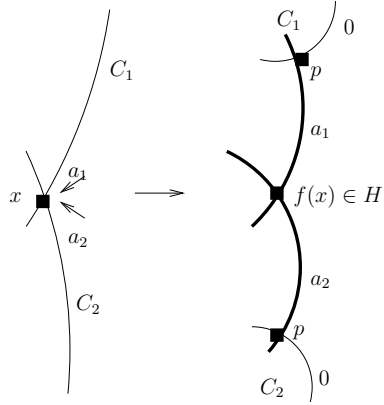}
\end {center}
\caption {The limit of the $\cs$ flow when a node maps to $H$.}
\label {nodal}
\end {figure}

For example, figure $\ref {nodal}$ shows the limit in the case of a node $x$ mapping to $H$ with contact orders $a_1$ and $a_2$ on the two components $C_1$ and $C_2$ transversal to $H$. The node is replaced by two rational components of degrees $a_1$ and $a_2$. These components are joined to the rest of the domain $C_1\cup C_2$ at nodes mapping to $p$. \\

We obtain the following algorithm for computing the limit $F$:
\begin {itemize} 
\item [(i)] We consider each irreducible component of the domain individually. We mark the nodes on each such component.
\item [(ii)] We leave unaltered the irreducible components mapping to $H$. 
\item [(iii)] The components which are transversal to $H$ are replaced in the limit by reducible maps. The reducible map has one back-bone component mapped to $p$. This component contains all markings which are not mapped to $H$. Moreover, rational tails are glued to the back-bone component at the points which map to $H$ according to the item below. The markings which map to $H$ are replaced by markings on the rational tails. 
\item [(iv)] In the limit, each isolated point $x$ of the domain curve which maps to $f(x)\in H$ with multiplicity $n$ is replaced by a rational tail glued at a node to the rest of the domain. The node is mapped to $p$. The image under $F$ of the rational tail is a curve in $\pr$ joining $p$ to $f(x)\in H$. The map $F$ is totally ramified over these two points with order $n$. If the point $x$ happens to be a section, we mark the point $F^{-1}(f(x))$ on the rational tail. 
\item [(v)]The map $F$ is obtained by gluing all maps in $(ii)$ and $(iii)$ along the markings we added in $(i)$ and then stabilizing.
\end {itemize}

\begin {corollary}\label{eqfam} For each stable map $f$, there is a family of stable maps $\eqref {family}$ over $\c$, whose fiber over $t\neq 0$ is the translated map $f^t$ and whose central fiber $F$ is obtained by the algorithm above.\end {corollary}

\begin {figure}
\begin {center}
\includegraphics*[scale=.7]{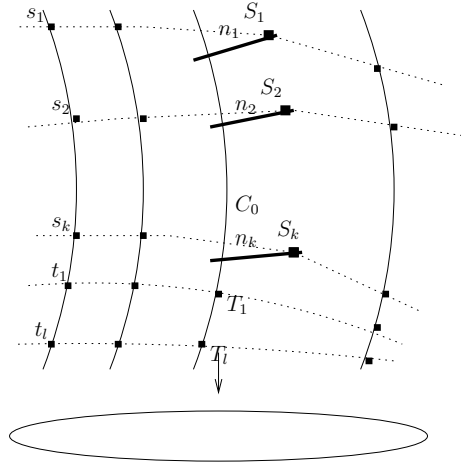}
\end {center}
\caption {A family as in corollary $\ref {eqfam}$.} 
\end {figure}

\subsection {Relation to the Gathmann stacks}\label{rel} We will proceed to identify the Bialynicki Birula cells of $\stack$. Recall that the fixed loci for the torus action on $\stack$ are indexed by decorated graphs $\Gamma$. We will identify the closed stacks $\overline {\mathcal F_{\Gamma}^{+}}$ in terms of images of fibered products of Kontsevich-Manin and Gathmann stacks under the tautological morphisms.

We will need the following versions of Gathmann's construction. 
\begin {itemize} 
\item [(i)] The substacks $\tgath$ of $\cgath$ parametrize maps with the additional condition that the components of $f$ are transversal to $H$. The maps in the open Gathmann stack $\ogath$ satisfy this condition by definition \cite{G}, hence: $$\ogath\hookrightarrow\tgath\hookrightarrow\cgath\hookrightarrow \stack.$$ 
\item [(ii-1)] For each map $f$ in $\tgath$, the dual graph $\Delta$ is obtained as follows:
\begin {itemize}
\item Vertices labeled by degrees correspond to the irreducible components of $f$. Vertices of degree $0$ satisfy the usual stability condition. 
\item The edges correspond to the nodes of $f$.  
\item Numbered legs correspond to the markings. The multiplicities $\alpha$ are distributed among the legs of $\Delta$. 
\item We write $\alpha_v$ for the ordered collection of multiplicities of the legs incoming to the vertex $v$ to which we adjoin $0$'s for all incoming edges (corresponding to the fact that the nodes of a map in $\tgath$ cannot be sent to $H$). The assignment of the multiplicities to the incoming flags is part of the datum of $\alpha_v$. 
\item The degree $d_v$ of the vertex $v$ is computed from the multiplicities: $|\alpha_v|=d_v$. 
\end {itemize}
For each graph $\Delta$ as above, we consider the stratum in $\tgath$ of maps whose dual graph is precisely $\Delta$. This is the image of the fibered product $\odgath$ of open Gathmann stack under the gluing maps: $$\odgath=\left(\prod_{v\in V(\Delta)} {\mathcal M}^{H}_{\alpha_v}(\pr, d_v)\right)^{E(\Gamma)}\to \tgath$$ The fibered product is computed along the evaluation maps on the corresponding moduli spaces at the markings determined by the edges of $\Delta$.

\item [(ii-2)] One also defines the stack $\widetilde {\mathcal M}_{\alpha}^{\Delta, H}(\pr, d)$ by taking the analogous fiber product of spaces $\widetilde {\mathcal M}^{H}_{\alpha_v}(\pr, d_v)$. Its image in $\tgath$ are the maps transversal to $H$ with domain type at least $\Delta$. We write $\overline {\mathcal M}^{\Delta}_{0,n}$ for the closure of the stratum of marked stable curves whose dual graph is the graph underlying $\Delta$ (forgetting the multiplicity labels). It follows that: \begin {equation}\label{fib} \widetilde {\mathcal M}_{\alpha}^{\Delta, H}(\pr, d)= \overline {\mathcal M}^{\Delta}_{0,n}\times_{\overline {\mathcal M}_{0,n}} \widetilde {\mathcal M}_{\alpha}^{H}(\pr, d). \end{equation} 

\item [(ii-3)]The similar fibered product of closed Gathmann spaces is defined as: \begin {equation}\label{cdgath}\cdgath=\left(\prod_{v\in V(\Delta)} {\mathcal {\overline M}}^{H}_{\alpha_v}(\pr, d_v)\right)^{E(\Gamma)}.\end {equation} We will see that $\cdgath$ is truly a compactification of $\odgath$. 

\item [(iii)] We will need to deal with unmarked smooth points of the domain mapping to $H$. This requires manipulations of a stack obtained from Gathmann's via the forgetful morphisms. We fix a collection of non-negative integers $\beta=(\beta_1, \ldots, \beta_n)$ and a collection of positive integers $\delta=(\delta_1, \ldots, \delta_m)$, satisfying the requirement $d=|\beta|+|\delta|$. We write ${\mathcal M}^{H}_{\beta, \delta}(\pr, d)$ for the image of the open Gathmann stack via the forgetful morphism: $${\mathcal M}^{H}_{\beta\cup \delta}(\pr, d)\hookrightarrow \overline {\mathcal M}_{n+m}(\pr, d)\to \stack.$$ The open stack ${\mathcal M}^{H}_{\beta, \delta}(\pr, d)$ parametrizes irreducible stable maps $f:C\to \pr$ with markings $x_1, \ldots, x_n$ such that: $$f^{\star} H= \sum \beta_i x_i + \sum \delta_j y_j,$$ for some distinct unmarked points of the domain $y_j$. 

\end {itemize}

With these preliminaries under our belt, we repackage the datum carried by each of the graphs $\Gamma$ indexing a fixed locus, into a fibered product $\mathcal X_{\Gamma}$ of Kontsevich-Manin and Gathmann spaces.  Precisely, we define: $$\mathcal X_{\Gamma}=\left(\prod_{v} \mathcal M^{H}_{\beta_v, \delta_v}(\pr, d_v)\times_{H}\prod_{w} \overline {\mathcal M}_{0, n(w)}(H, d_w)\right)^{E(\Gamma)}$$ The set of integers $\alpha_v$ defined in section $1.1$ is partitioned into two subsets $\beta_v\cup \delta_v$. $\delta_v$ collects the degrees of the incoming edges whose endpoint labeled $w$ is very unstable. The fibered product above is obtained as usual along the evaluation maps on the moduli spaces determined by the edges of $\Gamma$. 

A general point of $\mathcal X_{\Gamma}$ is obtained as follows.

\begin {itemize}  
\item For each vertex $w$ labeled $(H, d_w)$ we construct a stable map $f_w$ of degree $d_w$, with $n(w)$ markings and rational domain curve $C_w$. For the {\it unstable} vertices $w$ this construction should be interpreted as points mapping to $H$.
\item For each vertex $v$ labeled $p$, we construct a stable map $f_v$ with smooth domain $C_v$ and $n(v)$ marked points.
\item We join the domain curves $C_v$ and $C_w$ at a node each time there is an edge incident to both $v$ and $w$. Each edge $e$ which contains unstable $w$'s gives a special point of the domain. The special point should be a node mapping to $H$ if $w$ has two incoming edges, or a marking if $w$ has one incoming edge and an attached leg, or an unmarked point mapping to $H$ when $w$ is very unstable.
\item For each $v$ labeled $p$, the map $f_v$ has degree $d_v=|\beta_v|+|\delta_v|$ on the component $C_v$. Moreover, each incident edge $e$ corresponds to a point on $C_v$ which maps to $H$ and we require that the contact order of the map with $H$ at that point be $d_e$. \end {itemize} 
It is clear that by corollary $\ref {eqfam}$, the limit of the flow of the above map has dual graph $\Gamma$. 

The fibered product above is, modulo automorphisms, the Bialynicki Birula cell. We will carry out our discussion so that it only involves the stacks in $(i)$ and $(ii)$. To this end, we mark all the smooth points of the domain mapping to $H$. Combinatorially, this corresponds to eliminating the very unstable vertices in $\Gamma$. We let $\gamma$ be the graph obtained from $\Gamma$ by attaching one leg to each terminal very unstable vertex $w$. Then $\mathcal X_{\Gamma}$ is the image of the fibered product: $$\left(\mathcal X_{\gamma}=\right)\mathcal Y_{\Gamma}=\left(\prod_{v} \mathcal M^{H}_{\alpha_v}(\pr, d_v)\times_{H}\prod_{w} \overline {\mathcal M}_{0, n(w)}(H, d_w)\right)^{E(\Gamma)}$$ under the morphism: $$\overline {\mathcal M}_{0,n+\mathfrak u}(\pr, d)\to \stack$$ which forgets the markings corresponding to the $\mathfrak u$ newly added legs of $\gamma\to \Gamma$. We analogously define the companion stacks $\widetilde {\mathcal Y}_{\Gamma}$ and $\overline {\mathcal Y}_{\Gamma}$ (and their images $\widetilde {\mathcal X}_{\Gamma}$ and $\overline {\mathcal X}_{\Gamma}$): \begin {equation}\label {ytgamma} \widetilde {\mathcal Y}_{\Gamma}=\left(\prod_{v} \widetilde {\mathcal M}^{H}_{\alpha_v}(\pr, d_v)\times_{H}\prod_{w} \overline {\mathcal M}_{0, n(w)}(H, d_w) \right)^{E(\Gamma)}\end {equation} 
\begin {equation}\label {ygamma} \overline {\mathcal Y}_{\Gamma}=\left(\prod_{v} \overline {\mathcal M}^{H}_{\alpha_v}(\pr, d_v)\times_{H}\prod_{w} \overline {\mathcal M}_{0, n(w)}(H, d_w) \right)^{E(\Gamma)}.\end {equation}
There is a morphism $\overline {\mathcal Y}_{\Gamma} \to \overline {\mathcal M}_{0, n+\mathfrak u}(\pr, d)\to \stack$ obtained as compositions of: 
\begin {itemize}
\item gluing morphisms,
\item forgetful morphisms,
\item inclusions of Gathmann stacks $\overline {\mathcal M}^{H}_{\alpha}(\pr, d)\hookrightarrow \overline {\mathcal M}_{0, m}(\pr, d)$,
\item inclusions of Kontsevich-Manin stacks $\overline {\mathcal M}_{0,m}(H, d)\hookrightarrow \overline {\mathcal M}_{0, m}(\pr, d)$.
\end {itemize}

\begin {lemma}\label {dimsmooth}
The stack $\widetilde {\mathcal Y}_{\Gamma}$ is smooth. Its image in $\overline {\mathcal M}_{0, n+\mathfrak u}(\pr, d)$ has codimension $d+\mathfrak s-\mathfrak u$.  
\end {lemma}

\proof We observe that for any collection of weights $\alpha$, the evaluation morphism: \begin {equation}\label{sm}ev_1:\tgath\to H\end{equation} is smooth. First, the source is smooth. This is proved in \cite {G} for $\ogath$. To pass to the nodal locus, an argument identical to that of lemma 10 in \cite {FP} is required. As a consequence, there is a non-empty open set of the base over which the morphism is smooth. As $PGL(H)$ acts transitively on $H$, the claim follows. 

To prove the lemma, we follow an idea of \cite {KP}. We will induct on the number of vertices of the tree $\Gamma$, the case of one vertex being clear. We will look at the terminal vertices of $\Gamma$. 

Pick a terminal stable vertex $\mathfrak w$ labeled $(H, d_{\mathfrak w})$, if it exists. A new graph $\Gamma'$ is obtained by relabeling $\mathfrak w$ by $(H,0)$ and removing all its legs. Inductively, $\widetilde {\mathcal Y}_{\Gamma'}$ is smooth. It remains to observe that the morphism: $$\widetilde{\mathcal Y}_{\Gamma}\to{\widetilde{\mathcal Y}}_{\Gamma'}$$ is smooth, as it is obtained by base change from the smooth morphism: $$ev:\overline {\mathcal M}_{0, n(\mathfrak w)}(H, d_{\mathfrak w})\to H.$$ We can now assume all terminal vertices are either labeled $p$ or labeled $H$ but unstable. Removing all terminal $H$ labeled vertices from $\Gamma$, we obtain a new tree whose terminal vertices are all labeled $p$. Pick a terminal vertex $\mathfrak v$ in the new tree. It is connected to (at most) one vertex $\mathfrak w$. Assume $\mathfrak v$ was connected to the terminal vertices $\mathfrak w_1, \ldots, \mathfrak w_k$ in the old tree $\Gamma$. A new graph $\Gamma'$ is obtained from $\Gamma$ by removing all flags incident to $\mathfrak v, \mathfrak w_1, \ldots, \mathfrak w_k$ and replacing them by a leg attached at $\mathfrak w$.  The same argument as before applies. We base change $\widetilde {\mathcal Y}_{\Gamma'}$ by the smooth morphism $\eqref {sm}$. In our case, $\alpha=\alpha_{\mathfrak v}$ is the collection of degrees of the flags incoming to $\mathfrak v$. The evaluation is taken along the marking corresponding to the edge joining $\mathfrak v$ and $\mathfrak w$.
 
To compute the dimension of $\widetilde {\mathcal Y}_{\Gamma}$, we look at the contribution of each vertex $w$ labeled $(H, d_w)$, of each vertex $v$ labeled $p$, and we subtract the contribution of each edge $e$. Assuming all $H$ labeled vertices are stable, we obtain the following formula for the dimension of $\widetilde {\mathcal Y}_{\Gamma}$: $$\sum_{w} \left (rd_w+ (r-1) + n(w) - 3\right)+\sum_{v}\left((r+1) d_v + r + n(v) - 3- |\alpha_v|\right) - \sum_{e} (r-1) $$ $$=  r\left(\sum_{w} d_w+ \sum_{v} d_v\right)   +\left(\sum_{w} n(w)+\sum_{v}n(v)\right) + (r-4) W  + (r-3) V - (r-1) E$$ $$= rd + (n + 2E) + (-2E-W +r - 3)=\left((r+1)d+r+n-3\right)-d-W. $$ 

Thus the codimension of $\widetilde {\mathcal Y}_{\Gamma}$ in $\overline {\mathcal M}_{0,n}(\pr, d)$ equals $d+W$. The formula needs to be appended accordingly for the unstable vertices. The final answer for the codimension of $\widetilde {\mathcal Y}_{\Gamma}$ in $\overline {\mathcal M}_{0,n+\mathfrak u}(\pr, d)$ becomes $d+\mathfrak s - \mathfrak u$. 

We constructed open immersions $\mathcal Y_{\Gamma}\hookrightarrow \widetilde{\mathcal Y_{\Gamma}}\hookrightarrow \overline {\mathcal Y}_{\Gamma}.$ The following lemma clarifies the relationship between these spaces. 

\begin {lemma} \label{dense}
The image of $\mathcal Y_{\Gamma}$ is dense in ${\overline {\mathcal Y}}_{\Gamma}$. The stack ${\overline {\mathcal Y}}_{\Gamma}$ is reduced and irreducible.
\end {lemma}

\proof Using the above discussion, the only thing we need to show is that $\overline {\mathcal Y}_{\Gamma}$ is irreducible. We observe that the smooth stack $\mathcal Y_{\Gamma}$ is irreducible. Indeed, we can prove $\mathcal Y_{\Gamma}$ is connected by analyzing the $\cs$ action. Using corollary $\ref {eqfam}$ all maps in $\mathcal Y_{\Gamma}$ flow to one connected fixed locus which is the image of the connected stack: $$\prod_{v} \mathcal M_{0,n(v)}\times \prod_{w}\overline {\mathcal M}_{0, n(w)}(H, d_w).$$
To prove the irreducibility of $\overline {\mathcal Y}_{\Gamma}$ one uses the same arguments as in lemma $1.13$ in \cite {G}. By the above considerations, it is enough to show that any map $f$ in $\overline {\mathcal Y}_{\Gamma}$ can be deformed to a map with fewer nodes. Picking a map in $\overline {\mathcal  Y}_{\Gamma}$ is tantamount to picking maps  $f_v$ and $f_w$ in $\mathcal {\overline M}_{\alpha_v}^{H}(\pr, d_v)$ and in $\mathcal {\overline M}_{0, n(w)}(\pr, d_w)$ with compatible gluing data. For each vertex $v$, Gathmann constructed a deformation of $f_v$ over a smooth base curve such that the generic fiber has fewer nodes. We attach the rest of $f$ to the aforementioned deformation. To glue in the remaining components, we match the images of the markings by acting with automorphisms of $\pr$ which preserve $H$. The details are identical to those in \cite {G}. 

\begin {lemma}\label {weights}
There are $d+\mathfrak s - \mathfrak u$ negative weights on the normal bundle of $\mathcal F_{\Gamma}$. 
\end {lemma} 

{\bf Proof.} The arguments used to prove this lemma are well known (see for example \cite {GP} for a similar computation). In the computation below, we will repeatedly use the fact that the tangent space $T_{x}\pr$ has $\cs$ weights $D, \ldots, D$ for $x=p$ and weights $0, 0, \ldots, -D$ if $x\in H$. 

Recall the description of the stable maps in $\mathcal F_{\Gamma}$ which was given in the discussion following equation $\eqref {zeta}$. We let $(f, C, x_1, \ldots, x_n)$ be a generic stable map in $\mathcal F_{\Gamma}$ such that $C_v$ and $C_w$ are irreducible. We will compute the weights on the normal bundle at this generic point. These are the non-zero weights of the term $\mathcal T_f$ of the following exact sequence: $$ 0\to Ext^{0}(\Omega_{C}(\sum_{i}x_i), \mathcal O_{C}) \to H^{0}(C, f^{\star}T\pr)\to \mathcal T_{f} \to Ext^{1}(\Omega_{C}(\sum_{i}x_i), \mathcal O_{C})\to 0.$$ We will count the negative weights in the first, second and fourth term above. 

The first term gives the infinitesimal deformations of the marked domain. All contributions come from deformation of the components of type $C_e$. An explicit computation shows that the deformation space of such rational components with two special points, which need to be fixed by the deformation, is one dimensional with trivial weight. There is one exception in case the special points are not marked or nodes. This exceptional case corresponds to very unstable vertices. We obtain one negative weight for each such vertex, a total number of $\mathfrak u$.

Similarly, the fourth term corresponds to deformations of the marked domains. We are interested in the smoothings of nodes $x$ lying on two components $D_1$ and $D_2$. The deformation space is $T_xD_1\otimes T_x D_2$. The nodes lying on $C_e$ and $C_v$ give positive contributions. We obtain negative weights for nodes joining components $C_e$ and $C_w$ for stable $w$, and also for nodes lying on two components $C_{e_1}$ and $C_{e_2}$, which correspond to unstable $w$'s with two incoming edges. The number of such weights equals the number $F$ of edges whose vertex labeled $w$ is stable plus the number of unstable $w$'s with two incoming edges. 

The weights on the second term will be computed from the exact sequence $$0\to H^{0}(C, f^{\star} T\pr) \to \bigoplus_{v} H^{0}(C_v, f_v^{\star} T\pr) \bigoplus_{w} H^{0}(C_w, f_w^{\star} T\pr) \bigoplus_{e} H^{0}(C_e, f_e^{\star} T\pr)$$ $$\to\bigoplus_{\mathfrak f_v}  T_{\mathfrak f_v} \pr \bigoplus_{\mathfrak f_w} T_{\mathfrak f_w}\pr \to 0$$ Here $\mathfrak f_{v}$, $\mathfrak f_w$ are (some of the) flags of $\Gamma$ labeled by their initial vertices $v$ and $w$. They correspond to nodes of the domain mapping to $p$ and $H$, hence the terms $T_{\mathfrak f_v}\pr$ and $T_{\mathfrak f_w}\pr$ in the exact sequence above. The last term of the exact sequence above receives one negative contribution for each of the $F$ flags $\mathfrak f_w$. We obtain the following contributions to the negative weights of $H^{0}(C, f^{\star} T\pr)$ coming from the middle term. There are no negative contribution to $H^{0}(C_v,f_v^{\star}T\pr)=T_{p} \pr$. The Euler sequence: $$0\to \mathcal O\to \mathcal O_{\pr}(1)\otimes \c^{r+1}\to T\pr\to 0$$ and the arguments of \cite {GP} can be used to deal with the remaining two middle terms. Stable vertices labeled $w$ will contribute $d_w+1$ negative weights on $H^{0}(C_w,f_w^{\star}T\pr)$. Similarly there will be $d_e$ negative weights on $H^{0}(C_e, f_e^{\star}T\pr)$. We find that the number of negative weights of $H^{0}(C, f^{\star} T\pr)$ equals: $$\sum_{w} (d_w+1) + \sum_{e} d_e - F = d + \mathfrak s - F.$$ Thus the combined contributions of all the terms above to $\mathcal T_f$ is $d+\mathfrak s - \mathfrak u$. 
 
\begin {proposition}\label{gathmanncell} The closed cell $\overline {\mathcal F_{\Gamma}^{+}}$ is the stack theoretic image of the fibered product $\overline {\mathcal Y}_{\Gamma}$ of closed Gathmann and Konstevich-Manin spaces to $H$ under the tautological morphisms. Alternatively, it is the generically finite image of the stack $\overline {\mathcal X}_{\Gamma}$.\end {proposition}

\proof It is enough to show, by taking closures and using lemma $\ref {dense}$, that the stack theoretic image of $\widetilde {\mathcal Y}_{\Gamma}\to \stack$ is dense in $\mathcal F^{+}_{\Gamma}$. We observe that the geometric points of the image of $\widetilde {\mathcal Y}_{\Gamma}$ are contained in $\mathcal F_{\Gamma}^{+}$ because of corollary $\ref {eqfam}$. Moreover the dimensions match by lemmas $\ref {dimsmooth}$ and $\ref {weights}$. $\overline {\mathcal F_{\Gamma}^{+}}$ is reduced and irreducible because $\mathcal F_{\Gamma}$ clearly is, thanks to equation $\eqref {zeta}$. Same is true about $\tilde{\mathcal Y}_{\Gamma}$. These observations give our claim. The proof of proposition $\ref {filterability}$ shows that maps in $\overline {\mathcal Y}_{\Gamma}\setminus \widetilde {\mathcal Y}_{\Gamma}$ cannot flow to a map whose dual graph is $\Gamma$. As an afterthought, we obtain that the stack theoretic image of $\widetilde {\mathcal Y}_{\Gamma}$ equals $\mathcal F^{+}_{\Gamma}$.

\subsection {Filterability} In this subsection we define the partial ordering on the set of graphs indexing the fixed loci. This will be needed to establish the filterability condition $(c)$ of lemma $\ref {homologybasis}$. 

\begin {figure}
\begin {center}
\includegraphics*[scale=.6]{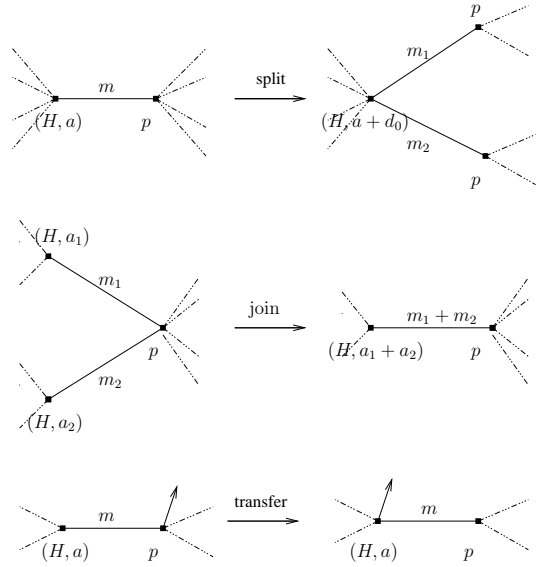}
\end {center}
\caption {Split, joins and transfers} 
\label {moves}
\end {figure}

\begin {figure}
\begin {center}
\includegraphics*[scale=.55]{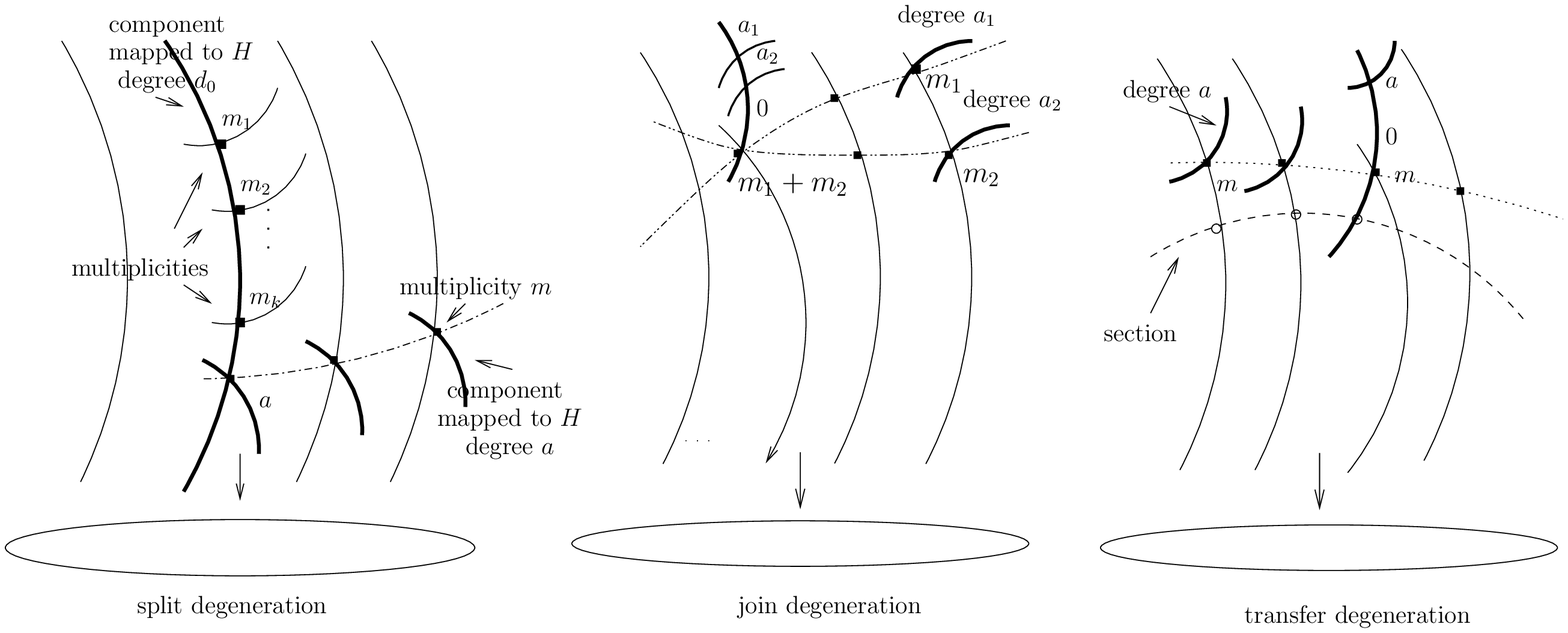}
\end {center}
\caption {Models for the combinatorial moves}
\label{localmoves}
\end {figure}

For any two decorated graphs $\Gamma$ and $\Gamma'$ indexing the fixed loci, we decree that $\Gamma\geq \Gamma'$ if there is a sequence of combinatorial surgeries called splits, joins and transfers changing the graph $\Gamma$ into $\Gamma'$. Each one of these moves is shown in figure $\ref {moves}$. Figure $\ref {localmoves}$ explains the intuition behind this ordering; we exhibit families of maps in a given Bialynicki-Birula cell degenerating to a boundary map which belongs to a different cell. The new cell should rank lower in our ordering. In figure $\ref {localmoves}$, the non-negative integers $a$ are degrees, and the positive integers $m$ are multiplicity orders with $H$. Components mapping to $H$ are represented by thick lines.

Explicitly, 
\begin {itemize}
\item A split move takes an edge of degree $m$ and cuts it into two (or several) edges with positive degrees $m_1$ and $m_2$. The vertex labeled $(H,a)$ is relabeled $(H, a+d_0)$ for some $d_0\geq 0$, while the vertex labeled $p$ is replaced by two vertices labeled $p$. The incoming edges and legs to the vertex $p$ are distributed between the newly created vertices. We require that $m=d_0+m_1+m_2$. The split move is obtained by degenerating a sequence of maps containing a point mapping to $H$ with multiplicity $m$. Such a degeneration is constructed in \cite {G}. The central fiber is a stable map in the boundary of the Gathmann space. There is an ``internal" component mapped to $H$ of degree $d_0$, to which other components are attached, having multiplicities $m_1, \ldots, m_k$ with $H$ at the nodes. The figure also shows an additional component mapped to $H$ with degree $a$ which is attached to the family.

\item The join move takes two edges of degrees $m_1$ and $m_2$ meeting in a vertex labeled $p$ and replaces them by a single edge whose degree is $m_1+m_2$, also collecting the two vertices labeled $H$, their degrees and all their incoming flags to a single vertex. Locally, the join move corresponds to a family of maps having two domain points mapping to $H$ with multiplicities $m_1$ and $m_2$ (there may be additional components mapping to $H$ with degrees $a_1$ and $a_2$ attached at these points). Letting the two points collapse, we obtain a boundary map with a point mapping to $H$ with multiplicity $m_1+m_2$. 

\item The transfer move can be applied to edges whose vertex labeled $p$ has an attached leg. We move the leg to the other end of the edge, labeled $H$. This move can be realized by a family of maps with one marking, and with domain points which map to $H$ with multiplicity $m$. In the limit, the marking and the point mapping to $H$ collapse. 

\end {itemize}

\begin {figure}
\begin {center}
\includegraphics*[scale=.7]{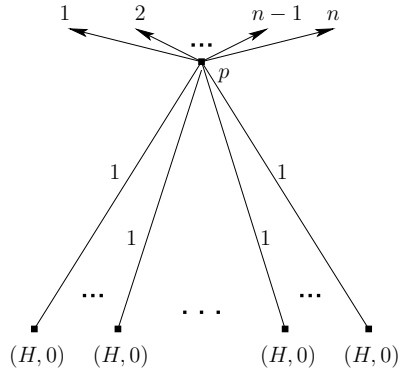}
\end {center}
\caption {The maximal graph.}
\label {maximal}
\end {figure}

To check that we have indeed defined a partial ordering we introduce the following length function:
\begin {eqnarray*} l(\Gamma)&=&\sum_{e} (e-1) \cdot\#\{\text {vertices labeled } (H,e) \}+ \# \{\text {vertices labeled } p\}+\\ &+& \#\{ \text {legs incident to }H\text{ labeled vertices}\}.\end {eqnarray*} The binary relation $"\geq"$ is indeed anti-symmetric since if $$\Gamma>\Gamma' \text { then } l(\Gamma)<l(\Gamma').$$ Moreover, it is clear that condition $(b)$ is satisfied; the unique maximal graph is shown in figure $\ref {maximal}$.

\subsection {The spanning cycles.} We will construct a family of cycles $\Xi$ satisfying the filterability condition $(c)$ of lemma $\ref {homologybasis}$. 

To begin with, we compare the cohomology and the Chow groups of the fixed loci, assuming that the filterability condition is satisfied. 

\begin {lemma}\label {isomorphic} 
The rational cohomology and rational Chow groups of $\mathcal F_{\Gamma}$ are isomorphic. The rational cohomology and rational Chow groups of $\stack$ are isomorphic. 
\end {lemma}

\proof We will use induction on $r$. There are two statements to be proved, we call them $A_r$ and $B_r$ respectively. It is proved in \cite {K} that $B_0$ is true. Lemma $\ref {homologybasis}$ shows that $A_r\implies B_r$. We conclude the proof by showing $B_{r-1}\implies A_r$. Indeed, the cohomology of $\mathcal F_{\Gamma}$ is computed using equation $\eqref{zeta}$ as: $$\left(\otimes_{v} H^{\star}(\overline {\mathcal M}_{0,n(v)})\otimes_{w} H^{\star}(\overline {\mathcal M}_{0, n(w)}(H, d_w))\right)^{A_{\Gamma}}.$$ It is remarkable that the same formula holds for the Chow groups. This follows from the Kunneth formula proved in theorem $2$ of \cite {K}. Our claim is established.

\begin {corollary}\label{hom} Let $X=G/P$ be any homogeneous space where $G$ is semisimple algebraic group and $P$ is a parabolic subgroup and $\beta\in A_1(X)$. Then parts $(i)$ and $(ii)$ of lemma $\ref {homologybasis}$ are true for $\fstack$.   
\end {corollary}

\proof We use a $\t$-action on $X$ with isolated fixed points. The fixed loci of the induced action on $\fstack$ are, up to a finite group action, products of the Deligne Mumford spaces $\overline {\mathcal M}_{0,n}$. The rational cohomology and Chow groups of the fixed loci are isomorphic. Using corollary $\ref {genatlas}$ and proposition $\ref {stratification}$ we obtain a Bialynicki-Birula decomposition on $\fstack$. We need to verify the conditions $(a)$ and $(b)$ required to apply lemma $\ref {homologybasis}$. We can check them on the closed points, hence we can pass to the coarse moduli schemes (considered in the sense of Vistoli \cite {V}). The two conditions are satisfied on general grounds on the projective irreducible \cite {KP} coarse moduli scheme $\fscheme$ of $\fstack$ as shown in \cite {B2}. We conclude observing that the image of $\mathcal F_i^{+}$ in $\fscheme$ is the corresponding Bialynicki-Birula cell $F_i^{+}$. In fact, one shows that $F_i^{+}$ is a coarse moduli scheme for $\mathcal F_i^{+}$. This is because $\mathcal F_i^{+}$ is reduced as $\mathcal F_i$ is reduced. 

The proof of lemma $\ref {isomorphic}$ also suggests the family $\Xi$. For each graph $\Gamma$, we will perform the following construction:
\begin {itemize}\item For each vertex $v$, pick a cycle class $\sigma_v$ on $\overline {\mathcal M}_{0, n(v)}$.
\item For each vertex $w$, pick a cycle class $\sigma_w$ on $\overline {\mathcal M}_{0,n(w)}(H, d_w)$. 
\item Assume that our choices define an $A_{\Gamma}$ invariant collection of classes.  
\end {itemize}

Henceforth, we will use explicit representatives for the above classes. Since $A^{\star}(\overline M_{0,n(v)})$ is generated by the boundary classes, we may assume:
\begin {itemize}\item $\sigma_v$ is the closed cycle $\overline {\mathcal M}^{\Delta_v}_{0,n(v)}$ of curves with dual graph (at least) $\Delta_v$. Here $\Delta_v$ is a stable graph with $n(v)$ labeled legs which are in one to one correspondence to the flags of $\Gamma$ incident to $v$. \end{itemize} In the following, $\xi$ will be any one of the cycles: \begin {equation}\label{cycle}\left[\prod_{v}  \overline {\mathcal M}^{\Delta_v}_{0,n(v)}\times \prod_{w} \sigma_w\big {/} A_{\Gamma}\right].\end {equation}

\begin {proposition}\label {filterability} The filterability condition $(c)$ is satisfied for the cycles $\xi$ defined above.
\end {proposition}
\proof By construction, it is clear that the cycles $\xi$ span the Chow groups of the fixed loci. 

We describe the closed points of $\xi^{+}$ informally. The dual graphs $\Delta_v$ determine the type of the domain curves. We consider maps with such domains which are transversal to $H$; points mapping to $H$ (with multiplicities determined by the edge degrees $d_e$ in $\Gamma$) are distributed on the irreducible components. Thus the closed points of $\xi^{+}$ give maps in smaller Gathmann spaces and the cycles $\sigma_w$.

Formally, we begin by adding one leg at each very unstable vertex of $\Gamma$, thus obtaining a graph $\gamma$ without very unstable vertices. Geometrically, this corresponds to marking {\it all} the smooth points on the domain which map to $H$. There is a forgetful morphism: \begin {equation}\label {for}\overline {\mathcal M}_{0, n+\mathfrak u}(\pr, d)\to \stack\end {equation} corresponding to the collapsing map $\gamma\to \Gamma$. 

A priori, the only decorations $\Delta_v$ carries are the labeled legs. The legs of $\Delta_v$ are in one-to-one correspondence to the incoming flags to $v$ in the graph $\Gamma$. However, we have seen in section $1.1$ that all flags of $\Gamma$ incident to $v$ carry the degrees $\alpha_v$. In this manner, we enrich the decorations of $\Delta_v$ using these degrees as ``multiplicities" associated to the legs (the half edges of $\Delta_v$ have multiplicity $0$). We denote by $\alpha(v)$ the datum of the collection of multiplicities $\alpha_v$ together with their distribution along the legs of $\Delta_v$. We can then form the fibered product $\widetilde {\mathcal M}_{\alpha(v)}^{\Delta_v, H}(\pr, d_v)$ as in section $\ref {rel}$ (ii-2). 

We let $\xi^{+}=\xi\times_{\mathcal F_{\Gamma}} \mathcal F_{\Gamma}^{+}$. We showed in proposition $\ref {gathmanncell}$ that there is a surjective morphism $\widetilde {\mathcal Y}_{\Gamma}\big/\text{Aut}_{\Gamma}\to \mathcal F_{\Gamma}^{+}$, inducing a surjective morphism $\xi\times_{\mathcal F_{\Gamma}}\widetilde{\mathcal Y}_{\Gamma}\big/\text{Aut}_{\Gamma}\to \xi^{+}$. We may need to endow these stacks with their reduced structure, but this suffices for the arguments in the Chow groups. By equations $\eqref {cycle}$, $\eqref{zeta}$, $\eqref {ytgamma}$ and $\eqref {fib}$, we derive that $\xi^{+}$ is the image of: $$\chi^{+}= \left[\left(\prod_{v} \widetilde {\mathcal M}_{\alpha(v)}^{\Delta_v, H}(\pr, d_v)\times_{H} \prod_{w} \sigma_w\right)^{E(\Gamma)} \big {/} \text{Aut}_{\Gamma}\right].$$ 
We similarly define: \begin {equation}\label{chi}\overline {\chi^{+}}= \left[\left(\prod_{v} \overline {\mathcal M}_{\alpha(v)}^{\Delta_v, H}(\pr, d_v)\times_{H} \prod_{w} \sigma_w\right)^{E(\Gamma)} \big {/} \text {Aut}_{\Gamma}\right].\end {equation} We let $\overline {\xi^{+}}$ be the its image under the forgetful map $\eqref{for}$. We obtain morphisms: $$\chi^{+}\hookrightarrow \overline {\chi^{+}}\to \overline {\mathcal Y}_{\Gamma}\big{/} \text{Aut}_{\Gamma}\to \overline {\mathcal F_{\Gamma}^{+}}.$$ The first one is an open immersion and by flatness of $\eqref {for}$ the same is true about the first inclusion below: $$\xi^{+}\hookrightarrow \overline {\xi^{+}}\to \overline {\mathcal F^{+}_{\Gamma}}.$$ We do not know that $\overline {\xi^{+}}$ is the closure of $\xi^+$ (we do not know $\overline {\xi^{+}}$ is irreducible). However, when formulating lemma $\ref {homologybasis}$ we were careful not to include this as a requirement in condition $(c)$.

Finally, we show that a map $f$ contained in the boundary $\overline {\xi^{+}}\setminus {\xi^{+}}$ flows to a fixed locus indexed by a graph $\Gamma'$ with $\Gamma'<\Gamma$. We first make a few reductions. Replacing $\Gamma$ by $\gamma$ and $\xi^{+}$ by $\chi^{+}$, we may assume $\Gamma$ has no very unstable vertices. We want to show that the graph of $F=\lim_{t\to 0} f^{t}$ is obtained from $\Gamma$ by a sequence of the combinatorial moves which we called joins, splits and transfers. 

\begin {figure}
\begin {center}
\includegraphics*[scale=.8]{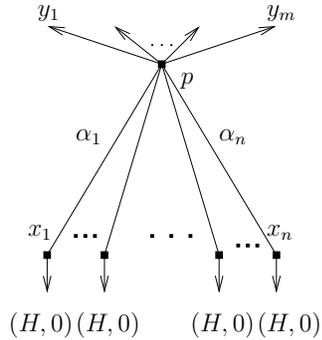}
\end {center}
\caption {The graph $\Gamma$ of a limit for a generic map in the Gathmann space.}
\label {gamma}
\end {figure}

The datum of a map $f$ is tantamount to giving maps $f_v$ and $f_w$ in the Gathmann spaces $\overline {\mathcal M}_{\alpha(v)}^{\Delta_v, H}(\pr, d_v)$ and the cycles $\sigma_w$ with compatible gluing conditions. As usual, the unstable vertices $w$ require special care as they only give points on the domain not actual maps. We have seen that the limit $F$ of $f^{t}$ is obtained from gluing the individual limits $F_v$ and $F_w=f_w$ (for stable $w$'s) of $f_v^{t}$ and $f_w^{t}$. The dual graphs are also obtained by gluing. Since to compute the limit we consider each vertex at a time, we may further assume that $\Gamma$ consists in one vertex labeled $v$ to which we attach legs and unstable vertices $w$. Thus, we may take $\Gamma$ to be the graph in figure $\ref {gamma}$.

Now recall that $\Delta_v$ encodes the domain type of the nodal map $f_v$. The markings of $f_v$ are distributed on the components of the domain and come with multiplicities encoded by the flags of $\Delta_v$. As we can treat the components individually, we may assume $\Delta_v$ has only one vertex. Moreover, the map $f_v$ has to be in the boundary of the Gathmann space $\overline {\mathcal M}^{H}_{\alpha(v)}(\pr, d_v)\setminus \widetilde {\mathcal M}^{H}_{\alpha(v)}(\pr, d_v).$ 

Changing notation slightly, we prove the following. We let $\alpha=(\alpha_1, \ldots, \alpha_n, 0, \ldots, 0)$ with $|\alpha|=d$ and $\alpha_i>0$. We will consider marked maps in the Gathmann space $\cgath$. For such a map $f$ with domain $C$ and markings $x_1, \ldots, x_n, y_1, \ldots, y_m$ we have $f^{!}H=\sum_{i} \alpha_i x_i.$ If $f$ were an element in $\tgath$, then its limit $F$ would have the dual graph $\Gamma$ shown in figure $\ref {gamma}$. This graph has one vertex $v$ labeled $p$, $n$ edges labeled $\alpha_1, \ldots, \alpha_n$ joining $v$ to unstable vertices $w$ labeled $(H,0)$. 
\begin {lemma}\label {comp}
Let $f$ be a map contained in the boundary of the Gathmann stack $\cgath\setminus \tgath$. Let $F$ be the limit of the flow of $f$. Then the dual graph $\Gamma_F$ of $F$ can be obtained from $\Gamma$ by splits, joins and transfers.\end {lemma}

\begin {figure}
\begin {center}
\includegraphics*[scale=.75]{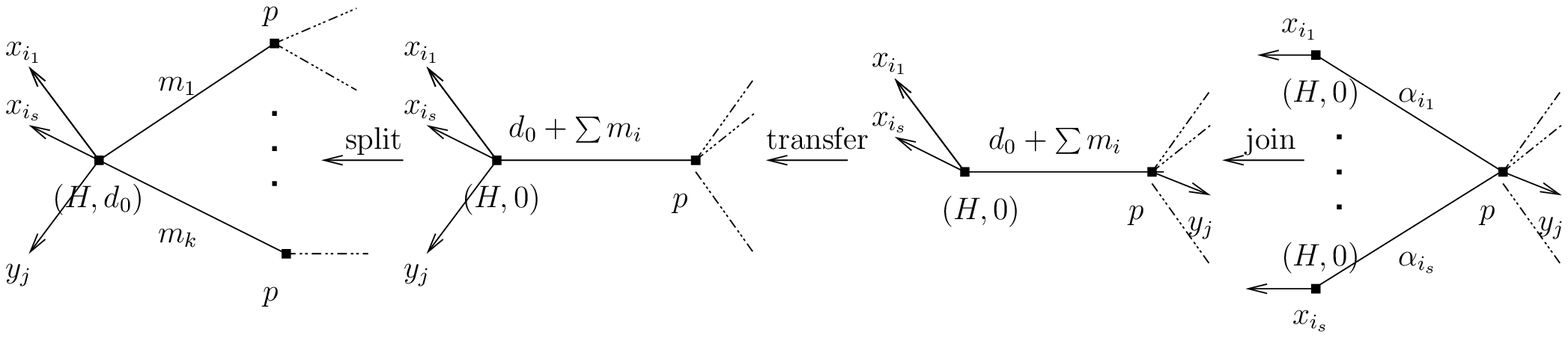}
\end {center}
\caption {The combinatorial moves comparing $\Gamma_F$ to $\Gamma$.}
\label {gammaf}
\end {figure}

\proof The map $f$ will have components which are not transversal to $H$ and which are responsible for the different dual graph. Let $C_0$ be a positive dimensional connected component of $f^{-1}(H)$ on which the map has total degree $d_0$, and let $C_1, \ldots, C_k$ be the irreducible components joined to $C_0$, having multiplicities $m_1, \ldots, m_k$ with $H$ at the nodes. Figure $\ref {boundary}$ shows an example of such a map. In any case, $C_0$ will contain some of the markings mapping to $H$, say $x_i$ for $i\in I$, and some of the remaining markings $y_j$ for $j\in J$. The contribution of the components $C_0\cup C_1 \ldots \cup C_k$ to the dual graph $\Gamma_F$, as computed by corollary $\ref {eqfam}$, is shown in the first graph of figure $\ref {gammaf}$. The figure also shows the moves (in reverse order) we apply to this portion of $\Gamma_F$ to obtain its corresponding contribution to $\Gamma$. The rest of the graph $\Gamma_F$ is attached to the portion shown there and is carried along when performing the combinatorial moves. Observe that existence of the join move is guaranteed by the equation $d_0+\sum m_i = \sum_{i\in I} \alpha_i$ which follows by considering intersection multiplicities of $f$ with $H$. The repeated application of this procedure will give a sequence of moves transforming $\Gamma$ into $\Gamma_{F}$. This completes the proof. 

\section {The tautology of the Chow classes}

In this section we tie the loose ends and prove the main result, theorem $\ref {main}$. Items (1)-(5) are contained in proposition $\ref {stratification}$, proposition $\ref {filterability}$, proposition $\ref {gathmanncell}$, lemma $\ref{weights}$ and lemma $\ref{isomorphic}$ respectively. Item $(6)$ is a consequence of the proof of proposition $\ref {filterability}$ and equations $\eqref{cdgath}$ and $\eqref{chi}$. 

\begin {lemma} \label {taut}

\begin {itemize}

\item [(i)] Let $i: H \to \pr$ be a hyperplane and let $i: \overline {\mathcal M}_{0,n}(H, d) \to \stack$ denote the induced map. The pushforward map $$i_{\star}:A_{\star}(\overline {\mathcal M}_{0,n}(H, d))\to A_{\star}(\stack)$$ preserves the tautological classes. 

\item [(ii)] For each $n$-multindex $\alpha$, the class of the Gathmann space $\left[\overline {\mathcal M}^{H}_{\alpha}(\pr, d)\right]$ is tautological. \end {itemize} \end {lemma}

\proof Consider the bundle $\mathcal B=R\pi_{\star} ev^{\star} \mathcal O_{\pr}(1)$ where $ev$ and $\pi$ are the universal evaluation and projection morphisms. This is a rank $d+1$ vector bundle on $\stack$. As usual, the equation of $H$ gives a section of $\mathcal B$ which vanishes precisely on $\overline {\mathcal M}_{0,n}(H, d)$. 

We claim that $$R^{\star}(\overline {\mathcal M}_{0,n}(H, d))\subset i^{\star}R^{\star}(\stack).$$ We need to check $i^{\star}R^{\star}(\stack)$ satisfies the two conditions of definition $1$. Invariance under pushforwards follows from standard manipulations of the projection formula. The second condition is immediate, as all classes $\alpha_H$ on $H$ are obtained by restrictions of classes $\alpha_{\pr}$ on $\pr$ and $$ev^{\star}\alpha_H=ev^{\star} i^{\star}\alpha_{\pr}=i^{\star}ev^{\star} \alpha_{\pr}.$$
Any tautological class $\alpha$ on $\overline {\mathcal M}_{0,n}(H, d)$ is the restriction of a tautological class $\beta$ on $\stack$. Therefore, $$i_{\star} \alpha = i_{\star} i^{\star} \beta = \beta \cdot c_{d+1}(\mathcal B).$$ 
It suffices to prove that $c_{d+1}(\mathcal B)$ is tautological. A computation identical to Mumford's \cite {M} using Grothendieck-Riemann-Roch shows: \begin {equation}\label {mumford} ch(\mathcal B)=\pi_{\star}\left(e^{ev^{\star} H} \cdot\left( \frac{c_1(\omega_{\pi})}{e^{c_1(\omega_{\pi})}-1}+i_{\star} P(\psi_{\star}, \psi_{\bullet})\right)\right).\end {equation} Here $P$ is a universal polynomial whose coefficients can be explicitly written down in terms of the Bernoulli numbers. The morphism $i$ is the codimension $2$ inclusion of the nodes of the fibers of the universal curve $\pi:\overline {\mathcal M}_{0, n\cup\{\diamond\}}(\pr, d)\to\stack$. Under the standard identifications, this can be expressed as union of images of fibered products : $$\overline {\mathcal M}_{0, S_1\cup \{\star\}}(\pr, d_1)\times_{\pr} \overline {\mathcal M}_{0, \{\star, \diamond, \bullet\}}(\pr, 0) \times_{\pr}\overline{\mathcal M}_{0, \{\bullet\}\cup S_2}(\pr, d_2)\stackrel{i}{\rightarrow} \overline {\mathcal M}_{0, {S_1\cup S_2\cup \{\diamond\}}}(\pr, d)$$ for all partitions $S_1\cup S_2=\{1, \ldots, n\}$ and $d_1+d_2=d$. The classes $\psi_{\star}$ and $\psi_{\bullet}$ of equation $\eqref {mumford}$ are the cotangent lines at the markings $\star$ and $\bullet$.

To prove out claim, we need to argue that $c_1(\omega_\pi)$ and $\psi$ are tautological. This follows from the results of \cite {divisors}, where it is shown that all codimension $1$ classes are tautological.
 
An argument may be required to justify the application of the Grothendieck-Riemann-Roch theorem in our stacky context. There are several ways to go about this. One can argue on the coarse moduli schemes using \cite {V}. Alternatively, lemma 2.1.1 in \cite {divisors} shows that if $r$ is large enough the locus of maps with automorphisms has codimension at least $d+2$. Its complement is a fine moduli scheme $\mathcal M^{\star}$ and we can apply GRR for the universal morphism over $\mathcal M^{\star}$. Since we are only interested in $c_{d+1}(\mathcal B)\in A^{d+1}(\stack)=A^{d+1}(\mathcal M^{\star})$, the formula $\eqref {mumford}$ holds up to codimension $d+1$. To deal with the small values of $r$, we pick $N$ large enough and use the inclusion $j:\stack \hookrightarrow \overline {\mathcal M}_{0, n}(\mathbb P^{r+N}, d)$. The class in question can be expressed as a pullback of a class we already know to be tautological: $$\left[\overline {\mathcal M}_{0,n}(H, d)\right]=j^{\star} \left[\overline {\mathcal M}_{0,n}(\mathcal H, d)\right]\in j^{\star} R^{d+1}(\overline {\mathcal M}_{0, n}(\mathbb P^{r+N}, d))=R^{d+1}(\stack).$$ Here $\mathcal H\hookrightarrow \mathbb P^{r+N}$ is a hyperplane which intersects $\pr\hookrightarrow \mathbb P^{r+N}$ along $H$. This proves the claim.

Part $(ii)$ is a consequence of equation $\eqref {gathmann}$, using induction on the multindex $\alpha$. The correction terms are pushforwards of classes on the boundary strata. These classes are either lower dimensional Gathmann spaces or Kontsevich-Manin spaces to $H$ which are tautological by induction and by part $(i)$ of the lemma respectively. This finishes the proof of the lemma.

Lemma $\ref{taut}$ gives the last claim $(7)$ of theorem $\ref{main}$ when the graph $\Gamma$ has only one vertex. The general case poses the following difficulty: we do not know that the stacks $\overline {\chi^{+}}$ obtained in proposition $\ref{filterability}$ are irreducible and of the dimension given by the naive count. The argument of lemma $\ref{dimsmooth}$ only gives this claim for the open substack $\widetilde {\chi^{+}}$ (defined just as $\overline {\chi^{+}}$ but by means of the open Gathmann substacks $\tgath$). Thus, we cannot immediately conclude that the Chow classes $\left[\overline{\chi^{+}}\right]$ are the {\it refined Gysin gluing} of the the fiber product factors.

To deal with this inconvenience, we argue as follows. It is shown in \cite {O2} that an additive system of generators for the tautological rings is given by (closures of) cycles of maps with fixed topological domain type and incidence conditions to various subspaces of $\pr$ at the markings and nodes. For general choices these substacks have the right dimension. We may assume $\sigma_w$ are represented by cycles of this form. Next, we use the arguments in section $2$ of \cite {G} to decrease all non-zero multiplicities $\alpha(v)$ to $1$, for each vertex $v$. This reduction is achieved by an equation of type $\eqref{gathmann}$. This fact is clear on the open stratum $\widetilde{\chi^{+}}$ where the relevant substacks behave as expected. Moreover, the boundary terms flow into fixed loci $\mathcal F_{\Gamma'}$ with $\Gamma'<\Gamma$. This was shown in proposition $\ref{filterability}$, and it was the technical assumption which made lemma $\ref{homologybasis}$ $(c)$ work. Finally, once all multiplicities carried by the $v$ vertices are either $0$ or $1$ (that is, if we only assume incidence conditions to $H$ at the markings or nodes of the domain curve indexed by $v$), the proof of lemma $\ref{taut}$ $(i)$ goes through as before.

\end {document}